
\tolerance=10000
\raggedbottom

\baselineskip=15pt
\parskip=1\jot

\def\sk{\vskip 3\jot}

\def\heading#1{\vskip3\jot{\noindent\bf #1}}
\def\label#1{{\noindent\it #1}}
\def\QED{\hbox{\rlap{$\sqcap$}$\sqcup$}}


\def\ref#1;#2;#3;#4;#5.{\item{[#1]} #2,#3,{\it #4},#5.}
\def\refinbook#1;#2;#3;#4;#5;#6.{\item{[#1]} #2, #3, #4, {\it #5},#6.} 
\def\refbook#1;#2;#3;#4.{\item{[#1]} #2,{\it #3},#4.}


\def\({\bigl(}
\def\){\bigr)}


\def\al{\alpha}
\def\be{\beta}

\def\ep{\varepsilon}
\def\ze{\zeta}
\def\et{\eta}

\def\ka{\kappa}
\def\la{\lambda}

\def\ch{\chi}

\def\Ga{\Gamma}
\def\De{\Delta}



\def\bfN{{\bf N}}

\def\bfR{{\bf R}}

\def\bfZ{{\bf Z}}

\def\abs#1{\big\vert#1\big\vert}
\def\norm#1{\vert#1\vert}

\def\oes{opposite-edge set}

{
\pageno=0
\nopagenumbers
\rightline{\tt low.fault.rate.arxiv.tex}
\vskip1in

\centerline{\bf Fault Tolerance in Cellular Automata at Low Fault Rates}
\vskip0.5in

\centerline{Mark McCann}
\centerline{\tt mark.mccann@asu.edu}
\centerline{GeoData Center for Geospatial Analysis and Computation}
\centerline{Arizona State University}
\centerline{Lattie F. Coor Hall} 
\centerline{P. O. Box 87530}
\centerline{Tempe, AZ 85287-5302}
\vskip0.5in

\centerline{Nicholas Pippenger}
\centerline{\tt njp@hmc.edu}
\centerline{Department of Mathematics}
\centerline{Harvey Mudd College}
\centerline{1250 Dartmouth Avenue}
\centerline{Claremont, CA 91711}
\vskip0.5in

\noindent{\bf Abstract:}
A commonly used model for fault-tolerant computation is that of cellular automata.
The essential difficulty of fault-tolerant computation is present in the special case of
simply remembering a bit in the presence of faults, and that is the case we treat in this paper.
The conceptually simplest mechanism for correcting errors in a cellular automaton
is to determine the next state of a cell by taking a majority vote among its neighbors
(including the cell itself, if necessary to break ties).
We are interested in which regular two-dimensional tessellations can tolererate faults using this mechanism,
when the fault rate is sufficiently low.
We consider both the traditional transient fault model (where faults occur independently in time and space) and a recently introduced combined fault model which
also includes manufacturing faults (which occur independently in space, but which affect cells for all
time).
We completely classify regular two-dimensional tessellations 
as to whether they can tolerate
combined transient and manufacturing faults, transient faults but not manufacturing faults, or not even transient faults.
\vfill\eject
}

\heading{1.  Introduction}

The results of this paper complement those of our previous paper, McCann and Pippenger [M3].
In that paper, we studied the fault-tolerance of certain cellular automata at high fault rates (close to 
$1/2$), and in this one we study the case of low fault rates (close to $0$).
In this introduction, we shall quickly describe the models we shall use; fuller discussion and motivation can be found in [M3].

In the theoretical study of fault-tolerant computation, the cellular automaton is one of the most frequently studied models.
This status comes about in the following way.
The most important theoretical model for computation, the Turing machine, is unsuitable for the 
study of fault-tolerant computation because it leaves large amounts of data on its tape unattended
for long periods of time while its head is elsewhere.
If each tape cell suffers a fault with fixed probability at each time step, the head (which can alter
the contents of only one cell per time step) has no hope of correcting these faults.
To evade this problem, it is natural to endow each tape cell with the ability to perform simple computations for the purpose of correcting local errors.
When this is done, the tape becomes a cellular automaton, and it is then natural to have the computation
performed by the cellular automaton in parallel, rather than serially by a moving head.
Deterministic cellular automata were introduced by Ulam [U] and von Neumann [N2].
Their probabilistic counterparts were introduced by Stavskaya and Pyatetski\u{\i}-Shapiro [S],
and were used to study fault-tolerant memory by Toom [T1, T2].

In the study of fault-tolerant computation by cellular automata, it is natural to separate the functions of fault-tolerant memory and computation.
In light of this separation, it is sufficient to consider the fault-tolerant cellular automata that remember a single bit, for they can then be combined with universal cellular automata for computation to obtain
universal fault-tolerant cellular automata. 
This idea is due to G\'{a}cs and Reif [G1].
A suitable two dimensional universal cellular automaton  with just two states per cell (the ``Game of Life'')
has been invented by Conway (see Gardner [G2], or Berlekamp, Guy and Conway [B]).
Thus we shall focus attention on {\it binary\/}  cellular automata, in which each cell has just two possible states.

For the correction of errors in binary cellular automata, majority voting is the most natural mechanism, and we 
will restrict attention in this paper to symmetric majority voting, in which the next state of each cell
is the result of a majority vote among the current states of the cells neighbors, including the current state of the cell itself if the number of neighbors is even.
This mechanism has many pleasant properties stemming from the fact that the majority function
is {\it monotone\/} (changing an input from $0$ to $1$ cannot change the output from $1$ to $0$) and  {\it self-dual\/} (complementing all the inputs complements the output).
Error-correction by majority voting was first analyzed by von Neumann [N1] in the context of Boolean circuits, and by Toom [T1, T2] for cellular automata.
The limitations of symmetric majority voting were discussed by Pippenger [P].

The cellular automata we consider are based on ``tessellations''.
A {\it tessellation\/} is an embedding of a graph into the sphere (if the graph is finite) or the plane
(if the graph is infinite).
The vertices and edges of the graph then divide the surface into which it is embedded into faces.
We focus on cellular automata based on {\it regular\/} tessellations, in which every vertex has the same degree $q$, and every face has the same degree $p$ (see Coxeter [C1]).
Thus the cellular automata we study have one cell at each vertex of the graph, with the neighbors of a cell being the cells located at adjacent vertices.
(Since the cells of a cellular automaton are in one-to-one correspondence with the vertices of its graph, we could if we wished use the terms ``cell'' and ``vertex'' interchangeably, or use one of these terms to the exclusion of the other.
We shall, however, usually use the term ``vertex'' when the context is graph-theoretic (that is, based on the edges and faces of the tessellation), and the term ``cell'' when the context is automata-theoretic
(that is, when the notions of state, fault or error are involved).)
Coxeter tessellations are denoted $\{p,q\}$, for $p\in\{3, 4, \ldots, \infty\}$ and $q\in\{2, 3, \ldots\}$.
(The preceding sentence brings out the need to distinguish the use of $\{p,q\}$ to name a tessellation from its use to name a set with two elements (or one if $\{p,q\}$).
This notation is well established for both purposes, however, and the meaning will always be clear
from context.)
The graph of $\{p,q\}$ has $q$ faces (which are $p$-gons) meeting at each vertex;
if $p=\infty$, it is an infinite $q$-regular tree (that is, an infinite tree in which every vertex has degree $q$).
These tessellations fall into three classes, as follows.
For $1/p + 1/q > 1/2$, $\{p,q\}$ is a finite tessellation of the sphere, corresponding to one of the Platonic solids (or, for $p < \infty$ and $q=2$, to two polygonal faces meeting at the equator).
For $1/p + 1/q = 1/2$, it is an infinite tessellation of the Euclidean plane (or, for $p=\infty$ and $q=2$, a tessellation of a line in the Euclidean plane).
Finally, for $1/p + 1/q < 1/2$, it is an infinite tessellation of the hyperbolic plane, with $q$ $p$-gons meeting at each vertex (or, for $p=\infty$ and $q\ge 3$, an infinite $q$-regular tree).
As we shall see below, only the hyperbolic tessellations have any hope of being fault tolerant with symmetric majority voting.
Although the tessellations of greatest interest for us are most naturally viewed as lying in the hyperbolic plane (where, according to hyperbolic geometry, all edges can have the same length and all faces the same area), hyperbolic geometry will play no role in our proofs, which are completely combinatorial in their methods.
The first study of cellular automata based on hyperbolic tessellations of the plane appears to be that of Wagner [W], who investigated ``modular'' computers in which a ``module'' was located at each vertex of $\{5,5\}$.

We come now to the specification of the fault models we shall consider.
We assume that the automaton is trying to remember a bit $x\in\{0,1\}$.
In this case, all the cells are initialized to state $x$, and any cell that at any subsequent time is in the complementary state $\overline{x}$ is said to be in {\it error\/} at that time.
The most natural model for faults in a binary cellular automaton is this: at each cell and each time step, the state of the cell is complemented independently with some probability $\epsilon>0$.
(This model is called the {\it purely probabilistic\/} fault model.)
For various reasons (which are discussed more fully in McCann and Pippenger [M3]), 
it is customary to work with a model which incorporates an adversarial as well as a probabilistic aspect.
At each cell and each time step a {\it fault\/} occurs independently with some probability 
$\epsilon>0$.
An adversary, who knows the value of $x$ and the locations of faults for all places and times, then
gains control of the state of each cell at each time at which a fault occurs for that cell.
The goal of the adversary to to maximize the probability that an error occurs at a particular cell at a particular time.
(This is called the {\it adversarial\/} fault model.
It is equivalent to the model used by Toom [T1, T2], where the role of the adversary is played by taking a supremum over a set of measures.)
For a binary automaton with a monotone transition rule (such as majority voting), the ``greedy strategy'' (create an error at any time and place at which a fault occurs) simultaneously maximizes the probability of error at all times and all places.
Furthermore, for a self-dual transition rule (such as majority voting), it is only necessary to consider
one of the two possible values of $x$ (since self-duality then allows us to deduce the behavior of the automaton for the other value).
We shall therefore always assume that $x=0$, so that a particular cell is in error at a particular time if and only if it is in state $1$ at that time.

The faults described in the preceding paragraph are {\it transient\/} faults, because if a fault occurs at a particular cell at a particular time, it does not affect the probability that a fault will occur at that cell at any future time.
McCann [M2] has introduced a fault model that is appropriate for studying {\it manufacturing\/} faults,
where each such fault occurs initially (during ``manufacturing'') and affects a cell for all time (that is, gives control of the state of the cell to the adversary for all time).
We shall study a {\it combined\/} fault model, in which each cell independently suffers a manufacturing fault with probability $\beta>0$, and suffers a transient fault at each time independently with probability 
$\alpha>0$.
We shall usually refer only to the single parameter 
$\epsilon = 1 - (1-\al)(1-\be) = \alpha + \beta - \alpha\beta$, 
which is the probability that a particular cell is controlled by the adversary at a particular time.
For either transient or combined faults, we shall say that a cellular automaton is {\it fault tolerant\/} if the 
probability that a particular cell is in error at a particular time is bounded by some quantity $\delta<1/2$
whenever $\epsilon>0$ is sufficiently small.
(In fact, in our positive results, we shall see that we can take $\delta = O(\epsilon)$ as $\epsilon\to 0$.)

Our results give a complete classification of the regular two-dimensional tessellations into three classes:
those that can tolerate combined faults (for a binary cellular automaton with symmetric majority voting), those that cannot tolerate combined faults but can tolerate transient faults, and those that cannot even tolerate transient faults.
For brevity, we shall speak of a tessellation tolerating faults of a given type, meaning that a binary cellular automaton based on that tessellation, and using symmetric majority voting, tolerates faults of that  type when the fault rate $\epsilon>0$ is sufficiently small.
This classification is summarized in the following table.

{\def\no{\cdot}\def\tro{\times}\def\bo{\otimes}
{\baselineskip=18pt
$$\matrix{
&&&&\ q&\rightarrow \cr
\cr
&            &2     &3      &4     &5     &6     &7     &8     &9     &\cdots \cr
&3          &\no &\no   &\no  &\no &\no &\tro  &\tro  &\bo  &\cdots  \cr
&4            &\no &\no &\no  &\tro  &\tro &\bo  &\bo  &\bo  &\cdots \cr
p\ \ &5            &\no &\no &\no  &\bo  &\bo &\bo  &\bo  &\bo  &\cdots \cr
\downarrow\ \ &6            &\no &\no &\no  &\bo  &\bo &\bo  &\bo  &\bo  &\cdots \cr
&\vdots \cr
&\infty       &\no &\tro &\tro  &\bo  &\bo &\bo  &\bo  &\bo  &\cdots \cr
\cr
}\eqno(1.1)$$
}
Here $\bo$ means ``tolerates combined faults'', $\tro$ means ``tolerates transient faults, but not combined faults'', and $\no$ means ``does not even tolerate transient faults''. }

Our results on transient faults are ultimately based on two established techniques:
for negative results, the notion of a ``self-sustaining island''
(implicit in the work of Toom [T2]), and for positive results, a theorem
(also due to Toom [T2]) that gives a sufficient condition for a cellular automaton to tolerate transient faults.
To apply Toom's theorem, however, we have had to develop a new technique, based on edge-colorings of tessellations that satisfy a certain invariance property:
all shortest paths from the origin to a given vertex have the same distribution of colors.
We call these edge-colorings ``shortest-path-invariant addressing schemes''.
(Schemes for describing paths in hyperbolic tessellations have been given for $\{5,5\}$ by Wagner [W] and for $\{4,5\}$ and $\{5,4\}$ by Marganstern [M1], but these schemes do not possess the shortest-path-invariance that we need for the application of Toom's theorem.)
Even this additional technique by itself, however, does not suffice for the tessellations $\{3,q\}$ with  $q\ge 7$.
In this case, the automaton does not appear to make enough progress in error correction during one time step to allow a direct application of Toom's theorem, though it does do so in two time steps.
We have thus been led to introduce the notion of the ``$\ka$-fold speed-up'' of a cellular automaton
(which accomplishes in one step what the original automaton accomplishes in $\ka$ steps), and to show that if the sped-up automaton tolerates transient faults, then so does the original automaton.

Our results on combined faults are based on entirely new methods.
Our negative results for combined faults use the notion of a ``pier-supported bridge'', analogous to the notion of a self-sustaining island that was used for transient faults.
Our positive results for combined faults cannot use Toom's theorem, whose proof depends on the transience of the faults in an essential way.
Instead, we prove an analogous theorem that gives a sufficient condition for a cellular automaton to tolerate combined faults.
To apply it, we use a ``balance of payments'' argument that is based not an edge-coloring (as was the case for Toom's theorem), but rather on an assignment of flows to edges.

In Section 2, we shall develop the combinatorial facts concerning regular hyperbolic tessellations that we shall need in subsequent sections, including the edge colorings that we shall use in our applications of Toom's theorem.
Our results on transient faults are presented in Section 3, and those on combined faults are presented in Section 4.
\sk

\heading{2. Regular Hyperbolic Tessellations}

All of our positive results apply to hyperbolic tessellations; spherical and Euclidean tessellations  are covered by our negative results: they do not even tolerate transient faults.
In this section, we collect various combinatorial facts about regular hyperbolic tessellations that we shall need in later sections.
The tessellation with faces of degree $p$ and vertices of degree $q$ is denoted $\{p,q\}$.
This tessellation is hyperbolic if and only if $1/p + 1/q < 1/2$.
This inequality is satisfied for $p=\infty$ and $q\ge 3$, for $p\ge 5$ and $q\ge 4$, for $p=4$ and 
$q\ge 5$, and for $p=3$ and $q\ge 7$.

Let $G = (V, E)$, with {\it vertices\/} $V$ and {\it edges\/} $E$, be the undirected graph of the tessellation $\{p,q\}$.
Two vertices $a$ and $b$ in $V$ are {\it adjacent}, or {\it neighbors\/} of each other, if they are joined by an edge $\{a,b\}$ in $E$.
Every vertex in $G$ has degree $q$.
The tessellation consists of an embedding of $G$ into the plane, in such a way that every 
{\it face\/} of the embedding has degree $p$.

Let us choose a vertex as the {\it origin\/} in $\{p,q\}$.
This choice partitions the vertices into sets that we shall call {\it generations\/} according to their distance (as measured by the number of edges in a shortest path) from the origin.
Generation $0$ comprises just the origin; for $g\ge 1$, generation $g$ comprises those vertices that are accessible through a path of length $g$, but no shorter path, from the origin.
Every edge either joins two vertices in consecutive generations, or joins two vertices in the same generation.
If vertex $a$ in generation $g\ge 0$ is joined by an edge with vertex $b$ in generation $g+1$,
we shall say that $a$ is a {\it parent\/} of $b$, or that $b$ is a {\it child\/} of $a$, and that the edge
is a {\it parent-child\/} edge.

The simplest regular hyperbolic tessellations are the tessellations $\{\infty,q\}$ with $q\ge 3$, which have no finite cycles, and thus are infinite trees.
The origin has no parents and $q$ children; every other vertex in generation $g\ge 1$ has one parent in generation $g-1$ and $q-1$ children in generation $g+1$.
Every edge is a parent-child edge.

The remaining tessellations $\{p,q\}$ are those with finite $p\ge 3$.
Their study breaks into two parts, according to whether $p$ is even or odd.
If $p\ge 4$ is even, every face, and therefore every cycle, contains an even number of edges,
so that the graph is bipartite.
Its vertices may therefore be colored with two colors in such a way that every edge joins vertices with different colors.
Since every vertex in generation $g\ge 1$ is joined by edges with vertices in the previous and subsequent generations,
this two-coloring must assign one color to vertices in even generations and the other to vertices in odd generations.
Thus no edge can join vertices in the same generation, so every edge is a parent-child edge.

On the other hand, if $p\ge 3$ is odd, there must be edges joining vertices in the same generation,
since otherwise coloring vertices according to the parity of their generations would yield a 
two-coloring of a graph containing odd cycles.
Edges joining vertices in the same generation $g\ge 1$ will be called either {\it sibling\/} edges
(if $p=3$, so that the vertices have a common parent in generation $g-1$), or {\it cousin\/} edges
(if $p = 2r+1\ge 5$, so that the vertices have a most recent common ancestor in 
generation $g-r\le g-2$).
(One should not read too much into this terminology:
vertices that are siblings have just one parent in common, not two,
and they also have a child in common!)
Vertices joined by a sibling (respectively, cousin) edge will be called {\it siblings\/} (respectively, 
{\it cousins\/}) of each other.

The simplest regular tessellations $\{p,q\}$ with even $p$ are those with $p=4$, which are hyperbolic for $q\ge 5$.
Consider two consecutive vertices $a$ and $b$ in generation $g\ge 1$, with $a$ to the left of $b$.
(We shall say ``to the left'' (respectively, ``to the right'') to mean counterclockwise (respectively, clockwise) in the cyclic order of a generation.
By the ``predecessor'' (respectively, ``successor'') of a vertex, we shall mean the next vertex to its left (respectively, right) in the cyclic order of a generation.
We shall say that two vertices are ``consecutive'' if one is the predecessor of the other.)
The vertices $a$ and $b$ lie on a common face of degree $4$, two of whose edges must join $a$ and $b$ to a common parent in generation $g-1$, and two of whose edges must join them to a common child in generation $g+1$.
This common child will therefore have two 
consecutive parents: it is the rightmost child of $a$ and the leftmost child of $b$.
Any vertex  that is neither a rightmost child nor a leftmost child will have just one parent.
Since a vertex can have at most two parents, it must have at least $q-2\ge 3$ children, and thus must have at least $3-2 = 1$ of these intermediate children.
We may therefore classify all vertices other than the origin into two classes:
{\it one-parent\/} vertices and {\it two-parent\/} vertices.
All the children of the origin are one-parent vertices; the leftmost and rightmost children of a vertex in generation $g\ge 1$ are two-parent vertices, and all its intermediate children are one-parent vertices.

The next simplest case is that of the tessellations $\{p,q\}$ with even $p = 2r\ge 6$.
These tessellations are hyperbolic when $q\ge 4$, but we shall only need to analyze those with 
$q\ge 5$, since for $q=4$ our results are negative ($\{p,4\}$ cannot tolerate even transient faults) and they can be proved by considering a single face.
Consider two consecutive vertices  $a$ and $b$ in generation $g\ge 1$, with $a$ to the left of $b$.
The vertices $a$ and $b$ lie on a common face of degree $p$, 
some of whose edges lie on a path $P$ between $a$ and $b$ though their most recent common ancestor,  and some of which lie on a path $Q$ through their earliest common descendant.
Since $Q$ has length at least two, $P$ has length at most $p-2$.
If $P$ has maximal length, $p-2$, then $Q$ has length two, so the rightmost child of $a$ will again be the leftmost child of $b$, and will again be a two-parent vertex.
If, however, $P$ has length less than $p-2$, then $Q$ will have length greater than two, so
the rightmost child of $a$ and the leftmost child of $b$ will be distinct one-parent vertices.

And as before, any children that are intermediate (neither leftmost nor rightmost) will be 
one-parent vertices.
Again, since a vertex can have at most two parents, it must have at least $q-2\ge 3$ children, and thus must have at least $3-2 = 1$ of these intermediate children.

\label{Lemma 2.1:}
For even $p\ge 4$ and $q\ge 5$, between a  two-parent vertex in a generation $g\ge1$ of 
$\{p,q\}$ and the next two-parent vertex to its right must lie at least $q-4\ge 1$ intermediate one-parent vertices.
In particular, at least one of the parents of a two-parent vertex (since they are consecutive) must be  a one-parent vertex.

\label{Proof:}
Let $a$ be a two-parent vertex, and let $b$ be the next two-parent vertex to its right.
Then $a$ must be the leftmost child of its right parent $c$.
Since $c$ has at most two parents, it must have at least $q-2\ge 3$ children.
The leftmost is $a$, and only the rightmost could be $b$.
Thus there are at least $q-2-2 = q-4\ge 1$
intermediate one-parent children between $a$ and $b$.
\QED

We turn now to the case of $\{p,q\}$ with odd $p\ge 3$.
The simplest of these tessellations are those with $p=3$, which are hyperbolic for $q\ge 7$.
Each pair of consecutive vertices in generation $g\ge 1$ is joined by a sibling edge, creating one face of degree three with their common parent in generation $g-1$, and another such face with their common child in generation $g+1$.
Each vertex in generation $g\ge 1$ has at most two parents and exactly two siblings, and thus has at least $q-2-2 = q-4 \ge 3$ children.
Only its leftmost and rightmost children are two-parent vertices, so there are at least $3-2 = 1$ intermediate one-parent children between $a$ and $b$.

If from $\{3,q\}$ we delete all sibling edges, we obtain another tessellation, which we shall denote 
$\{3,q\}'$, that is ``almost regular''.
In $\{3,q\}'$, every face has degree four (since it was created by deleting a sibling edge that was formerly incident with two faces of degree three).
The origin of $\{3,q\}'$ has degree $q$, but every other vertex has lost its two former siblings as neighbors, and thus now has degree $q-2$.
The tessellation $\{3,q\}'$ can therefore be obtained by partitioning the edges and vertices other than the origin of $\{4,q-2\}$ into $q-2$ isomorphic ``sectors'', then arranging $q$ of these sectors around the origin of $\{3,q\}'$.

\label{Lemma 2.2:}
For $q\ge 7$, between a  two-parent vertex in a generation $g\ge1$ of 
$\{3,q\}$ and the next two-parent vertex to its right must lie at least $q-6\ge 1$ intermediate one-parent vertices.
In particular, if $q\ge 8$, at least one of the siblings of a one-parent vertex must be another 
one-parent vertex.

\label{Proof:}
Let $a$ be a two-parent vertex, and let $b$ be the next two-parent vertex to its right.
Then $a$ must be the leftmost child of its right parent $c$.
Since $c$ has at most two parents, and exactly two siblings, it must have at least $q-2-2=q-4\ge 3$ children.
The leftmost is $a$, and only the rightmost could be $b$.
Thus there are at least $q-4-2 = q-6\ge 1$
intermediate one-parent children between $a$ and $b$.
\QED

Finally, we come to the last and most intricate case, that of $\{p,q\}$ with odd $p=2r+1\ge 5$.
These tessellations are hyperbolic when $q\ge 4$, but we shall only need to analyze those with 
$q\ge 5$, since for $q=4$ our results are negative ($\{p,4\}$ cannot tolerate even transient faults) and they can be proved by considering a single face.

Consider two consecutive vertices  $a$ and $b$ in generation $g\ge 1$, with $a$ to the left of $b$.
There are three possible cases.
Firstly, $a$ and $b$ may be joined by a path of length $p-1 = 2r$ through a most recent common ancestor in generation $g-r$.
In this case, $a$ and $b$ are joined by a cousin edge that completes a face of degree $p$,
and the rightmost child of $a$ is distinct from the leftmost child of $b$, with both of these children being one-parent vertices.
Second, $a$ and $b$ may be joined by a path of length $p-2 = 2r-1$ through a cousin edge joining two vertices in generation $g-r+1$.
In this case, $a$ and $b$ are not joined by a cousin edge, but the rightmost child of $a$ is also the leftmost child of $b$, completing a face of degree $p$, and this common child is a two-parent vertex.
Finally, if length of the shortest path from $a$ to $b$ through one or more vertices in earlier generations is at most $p-3 = 2r-2$, then $a$ and $b$ are not joined by a cousin edge,
and the rightmost child of $a$ is distinct from the leftmost child of $b$, with both of these children being one-parent vertices.

\label{Lemma 2.3:}
In $\{p,q\}$ with $p\ge 5$ and $q\ge 5$, a vertex can have at most two neighbors that are either its parents or its cousins.
In particular, each vertex must have at least three children.

\label{Proof:}
Suppose, to obtain a contradiction, that $a$ is a vertex with more than two neighbors that are parents or cousins, and that $a$ is in  the earliest generation $g\ge 1$ in which such a vertex can be found.
Then any parents of $a$ must have at least three children.
Suppose first that $a$ has two parents $b$ and $c$, with $b$ to the left of $c$.
Then $a$ is the rightmost child of its left parent $b$.
The predecessor $d$ of $a$ is therefore the second-rightmost child of $b$.
Thus $a$ cannot be joined by a cousin edge to its predecessor $d$, since that would create a face
with vertices $a$, $b$ and $d$ of degree three, contradicting the hypothesis that $p\ge 5$.
Similarly, $a$ cannot be joined by a cousin edge to its successor, since this would create a face of degree three with its right parent.
Thus, a two-parent vertex cannot have any cousins.

Now suppose $a$ has just one parent $b$.
By the reasoning above, $a$ can be joined by a cousin edge with its predecessor only if $a$ is the leftmost child of its parent $b$, and $a$ can be joined by a cousin edge with its successor only if  $a$ is the rightmost child of its parent $b$.
Since $b$ has at least three children, these two possibilities are mutually exclusive,
so a one-parent vertex can have at most one cousin.
\QED

We may thus classify one-parent vertices into two classes:
those incident with exactly one cousin edge will be called {\it cousin\/} vertices,
and those not incident with a cousin edge will be called {\it non-cousin\/} one-parent vertices.
By Lemma 2.3, every vertex other than the origin is a cousin vertex (and thus also a one-parent vertex), a non-cousin one-parent vertex, or a 
two-parent vertex.

\label{Lemma 2.4:}
For $p\ge 5$ and $q\ge 5$, a one-parent vertex in $\{p,q\}$ has at most one two-parent child.

\label{Proof:}
Suppose, to obtain a contradiction, that $a$ and $b$ are both two-parent children of the 
one-parent vertex $c$.
Suppose, without loss of generality, that $a$ is the leftmost child of $c$, and that $b$ is the rightmost child of $c$.
Let $d$ be the left parent of $a$; then $d$ is consecutive with $c$, which is the right parent of $a$.
Let $e$ be the right parent of $b$; then $e$ is consecutive with $c$, which is the left parent of $b$.
Let $f$ be the parent of $c$.
Neither $d$ nor $e$ can be a child of $f$, since doing so would create a face of degree four
($a$, $c$, $d$ and $f$, or $b$, $c$, $e$ and $f$),
contradicting the hypothesis that $p\ge 5$.
Thus $f$, which by Lemma 2.3 has at least three children, must have these children strictly between $d$ and $e$, contradicting the fact that $d$, $c$ and $e$ are consecutive, so that $c$ is the only vertex strictly between $d$ and $e$.
\QED

\label{Lemma 2.5:}
For $p\ge 5$ and $q\ge 5$, both parents of a two-parent vertex in $\{p,q\}$ are one-parent vertices.
Indeed, they are non-cousin one-parent vertices.

\label{Proof:}
Let $a$ be a two-parent vertex with left parent $b$ and right parent $c$.
Then $b$ and $c$ are consecutive.

Suppose first, to obtain a contradiction, that at least one of $b$ and $c$ is a two-parent vertex.
Suppose, without loss of generality, that $b$ is a two-parent vertex.
Let $d$ be the right parent of $b$.
By Lemma 2.3, $d$ has at least three children.
The leftmost is $b$, so the second-leftmost must be $c$ (since $c$ is the successor of $b$).
Then $a$, $b$, $c$ and $d$ form a face of degree four, contradicting the hypothesis that 
$q\ge 5$.
This contradiction shows that neither $b$ nor $c$ can be a two-parent vertex.

Suppose then, to obtain a contradiction, that at least one of $b$ and $c$ is a cousin vertex.
Suppose, without loss of generality, that $b$ has a cousin $e$.
The vertex $e$ must be the predecessor of $b$ (since $c$ is the successor of $b$, so if it were a cousin of $b$, then $a$, $b$ and $c$ would form a face of degree three, contradicting the hypothesis that $q\ge 5$).
Let $d$ be the parent of $b$.
By Lemma 2.3, $d$ has at least three children.
The leftmost must be $b$ (since if $e$ were a child of $d$, the vertices $b$, $d$ and $e$ would form a face of degree three, contradicting the hypothesis that $q\ge 5$).
Thus has a least two children to the right of $b$.
The second-leftmost must be $c$ (since $c$ is the successor of $b$).
Then $a$, $b$, $c$ and $d$ form a face of degree four, contradicting the hypothesis that 
$q\ge 5$.
This contradiction shows that neither $b$ nor $c$ can be a cousin vertex.
Thus both parents of a two-parent vertex must be non-cousin one-parent vertices.
\QED

\label{Lemma 2.6:}
Of two vertices  that are cousins of each other, at least one has a parent that is a non-cousin one-parent vertex.

\label{Proof:}
Suppose that vertices $a$ and $b$, with $a$ to the left of $b$, are cousins of each other.
Since $a$ and $b$ are cousins, we must have $q\ge 5$.
Let $c$ and $d$ be the parents of $a$ and $b$, respectively.
Since $a$ and $b$ are consecutive, so are $c$ and $d$, with $c$ to the left of $d$.
Suppose, to obtain a contradiction, that $c$ and $d$ are each either a two-parent vertex or a cousin vertex.

Consider first the case that both $c$ and $d$ are two-parent vertices.
Let $e$ be the right parent of $c$.
By Lemma 2.3,  $e$ has at least three children.
The leftmost is $c$, and only the rightmost could be $d$.
Thus there are at least $3-2=1$ children between $c$ and $d$, contradicting the fact that $c$ and $d$ are consecutive.

Next consider the case that one of $c$ and $d$ is a two-parent vertex, while the other is a cousin vertex.
Assume, without loss of generality, that $c$ is a two-parent vertex, while $d$ is a cousin vertex.
Let $f$ be the cousin of $d$.
Then $f$ must be the successor of $d$ (since if it were the predecessor $c$ of $d$,
then $a$, $b$, $c$ and $d$ would form a face of degree four, contradicting the fact that $q\ge 5$).
Let $e$ be the right parent of $c$.
By Lemma 2.3, $e$ has at least three children.
The leftmost is $c$, so the next two vertices to the right of $c$ (namely $d$ and $f$) must also 
be children of $e$.
So $d$, $e$ and $f$ form a face of degree three, contradicting the fact that $q\ge 5$.

Finally consider the case that $c$ and $d$ are both cousin vertices.
They cannot be cousins of each other, since then $a$, $b$, $c$ and $d$ would form a face of degree four, contradicting the fact that $q\ge 5$.
Thus the cousin $g$ of $c$ must be the predecessor of $c$, and the cousin $f$ of $d$ must be the successor of $d$.
Let $e$ be the parent of $c$.
By Lemma 2.3, $e$ has at least three children.
The leftmost must be $c$ (since if $g$ were a child of $e$, the vertices $c$, $e$ and $g$ would form a face of degree three, contradicting the fact that $q\ge 5$).
Thus the next two vertices to the right of $c$ (namely $d$ and $f$) must be children of $e$.
So $d$, $e$ and $f$ form a face of degree three, contradicting the fact that $q\ge 5$.
\QED

We turn now to the edge colorings that we shall use in our applications of Toom's theorem.
These will only be needed for $\{4,q\}$ with $q\ge 5$ and $\{3,q\}$ with $q\ge 7$, so we shall restrict our attention to $p\in \{3,4\}$.

An {\it addressing scheme\/} for a graph $G = (V,E)$ with $l$ colors is a function
$\ch : E\to \{0,\ldots, l-1\}$ that assigns a ``color''  $\ch(e)$ to each edge $e$ in $E$, and that assigns distinct colors to the edges in $E$ incident with a common vertex.
An addressing scheme for $G$ can be used to ``navigate'' in $G$:
a path from one vertex to another can be described by giving the sequence of colors of the edges traversed by the path.

For our results, we shall only need to deal with paths starting from a fixed origin, and only with shortest paths from the origin.
We shall, therefore, need colors only for edges that lie on shortest paths from the origin.
Let a vertex of $G$ be distinguished as the origin.
Let $E'\subseteq E$ be the set of edges lying on shortest paths from the origin to other vertices.
A {\it partial addressing scheme\/} for $G$ with $l$ colors is a function
$\ch : E'\to \{0,\ldots, l-1\}$ that assigns a ``color''  $\ch(e)$ to each edge $e$ in $E'$, and that assigns distinct colors to the edges in $E'$ incident with a common vertex.
(If $p=4$, then $E' = E$, so a partial addressing scheme is simply an addressing scheme;
if $p=3$, then $E'$ is obtained from $E$ by deleting all sibling edges, so only parent-child edges receive colors.)

To apply Toom's theorem, we shall need partial addressing schemes that satisfy an additional constraint.
A partial addressing scheme $\ch$ is {\it shortest-path-invariant\/} if, for any vertex $a$ in $V$, any two shortest paths $P$ and $Q$ from the origin to $a$, and any color $k\in \{0, \ldots, l-1\}$,
the number of edges with color $k$ in $P$ is the same as the number of edges with color $k$ in $Q$.
This common number will be denoted $\norm{a}_k$.
If the number of edges on a shortest path from the origin to $a$ (that is, the number of the generation containing $a$) is denoted $\norm{a}$, then we have 
$\norm{a} = \norm{a}_0 + \cdots + \norm{a}_{l-1}$, since each edge on a shortest path from the origin to $a$ has one of the $l$ colors.

\label{Lemma 2.7:}
The graph of thetessellation $\{\infty, q\}$ with $q\ge 3$ has a shortest-path-invariant addressing scheme using $q$ colors.

\label{Proof:}
Select an origin  and color the $q$ edges incident with the origin with the $q$ colors.
Every vertex in generation 1 now has just one incident edge (the edge joining it to its parent) colored,
and the remaining $q-1$ edges (the edges joining it to its children) can be colored with the remaining $q-1$ colors.
We may continue in this way, coloring the remaining edges of vertices in generations 
$2, 3, \ldots$ in succession, with the $q-1$ edges joining a given vertex to its children receiving the $q-1$ colors distinct from that of the edge joining it to its parent.
Since there is a unique path between any two vertices in a tree, the shortest-path-invariance
is satisfied trivially.
\QED

\label{Lemma 2.8:}
Let $G$ be the graph of a tessellation $\{4,q\}$, and
let a vertex of $G$ be distinguished as the origin.
Then a coloring of the edges of $G$ that assigns distinct colors to the edges incident with 
every vertex, and 
assigns the same color to opposite edges of every face,
is a shortest-path-invariant addressing scheme.

\label{Proof:}
Suppose, to obtain a contradiction, that $P$ and $Q$ are two shortest paths from the origin to a vertex $a$ in generation $g$, and that $P$ and $Q$ have different color distributions.
We may assume that $g\ge 2$ is the earliest generation for which shortest paths with different color distributions exist.

The paths $P$ and $Q$ must pass through distinct vertices in generation $g-1$, since otherwise
the common vertex $b$ of $P$ and $Q$ in generation $g-1$ would lie in an earlier generation and the paths $P$ and $Q$, truncated at $b$, would have different color distributions
(since the color distributions of the truncated paths are the same as those of $P$ and $Q$, but each with one fewer occurrence of the color $\ch(\{a,b\})$).

Thus the vertex $a$ must have two parents.
Let $b$ and $c$ be the left and right parents of $a$, respectively,
and suppose that $P$ passes through $b$ and $Q$ passes through $c$.
Since the parents $b$ and $c$ of $a$ are consecutive,
the vertices $a$, $b$ and $c$ lie on a face, and the fourth vertex of this face must be a common parent of both $b$ and $c$ in generation $g-2$.
Let $d$ be this common parent, and let $R$ be a shortest path from the origin to $d$.
Since $b$ lies in generation $g-1$, the path $P$, truncated at $b$, must have the same color distribution as the path $R$ extended by $\{b,d\}$.
Similarly, the path $Q$, truncated at $c$, must have the same color distribution as the path $R$ extended by $\{c,d\}$.

Thus $P$ has the same color distribution as that of the path $R$ extended by $\{b,d\}$ and 
$\{a,b\}$, and $Q$ has the same color distribution as that of the path $R$ extended by $\{c,d\}$ and  $\{a,c\}$.
But $\ch(\{a,b\}) = \ch(\{c,d\})$ and $\ch(\{a,c\}) = \ch(\{b,d\})$, so $P$ and $Q$ have the same color distribution.
This contradiction completes the proof.
\QED

\label{Lemma 2.9:}
The tessellation $\{4,q\}$ with $q\ge 5$ has a shortest-path-invariant addressing scheme using 
$q$ colors.

\label{Proof:}
Color an edge incident with the origin $0$, then color the remaining edges incident with the origin
$1, 2, \ldots, q-1$ in counterclockwise order.
This action colors one edge incident with each vertex in generation $1$.
For each vertex in generation $1$, color its remaining incident  edges with the remaining $q-1$ colors successively (modulo $q$) in clockwise order.
This action colors either one or two edges incident with each vertex in generation $2$, and if two edges
incident with a vertex have been colored, they will have successive colors (modulo $q$) in counterclockwise order.
For each vertex in generation $2$, color its remaining incident edges with the remaining $q-1$ or $q-2$
colors successively (modulo $q$) in counterclockwise order.
Continue in this way for generations $3, 4, \ldots$ as follows.
After all edges incident with vertices in generation $g$ have been colored, either one or two edges
incident with each vertex in generation $g+1$ will have been  colored, and if two have been colored,
they will have successive colors (modulo $q$) in counterclockwise order if $g+1$ is even, and in clockwise order if $g+1$ is odd.
It will then be possible, for each vertex in generation $g+1$, to color its remaining incident edges
with the remaining $q-1$ or $q-2$ colors successively (modulo $q$) in counterclockwise order if $g+1$ is even, and in clockwise order if $g+1$ is odd.
This procedure assigns distinct colors to the edges incident with a given vertex, and
(because of the alternation between ``clockwise'' and ``counterclockwise'' in consecutive generations) it assigns the same color to opposite edges of any face.
Thus by Lemma 2.8, it yields a shortest-path-invariant addressing scheme.
\QED

\label{Lemma 2.10:}
The tessellation $\{3,q\}$ with $q\ge 7$ has a shortest-path-invariant partial addressing scheme using $q$ colors.

\label{Proof:}
The deletion of all sibling edges from $\{3,q\}$ yields a tessellation which we shall denote
$\{3,q\}'$, in which the origin has degree $q$,
every vertex other than the origin has degree $q-2$ and  every face has degree  $4$.
The edges deleted do not lie on any shortest paths from the origin to any vertex, and thus it will suffice to assign colors to the parent-child edges in $\{3,q\}$, which are the edges of $\{3,q\}'$.

To color $\{3,q\}'$, we begin by coloring $\{4,q-2\}$ according to Lemma 2.9.
This colored $\{4,q-2\}$ may be embedded in the plane in a way that is symmetric under the dihedral group of order $2(q-2)$ of rotations and reflections fixing the origin, so that it is 
divided by $q-2$ ``spokes'' into $q-2$ ``sectors''.
The {\it spoke\/} $s$ (where $0\le s\le q-3$) will be a ray from the origin containing the edge from the origin 
colored $s$.
The origin will not belong to any spoke or sector; the edges and all other  vertices will be partitioned
into sectors.
If $q$ is even, a spoke will contain one edge directed out of each generation,
and one vertex in each generation other than the origin. 
If $q$ is odd, a spoke will 
contain one edge directed out of each generation congruent to $0$ modulo $3$,
and one vertex in each generation congruent to $0$ or $1$ modulo $3$.
The {\it sector\/} $s$ will comprise all edges and  vertices other than the origin between spoke $s$  (inclusive) and spoke $s+1$ (modulo $q-2$, exclusive).

The tessellation $\{3,q\}'$ can be obtained from $\{4,q-2\}$ by arranging $q$ sectors of $\{4,q-2\}$
(compressed angularly) around the origin, with new spokes which we shall number $q-2$ and $q-1$.
Thus $\{3,q\}'$ (apart from the origin) is partitioned $q$ sectors, which are numbered $0$ through $q-1$,
with sector $s$ (where $0\le s\le q-1$) comprising all edges and vertices other than the origin between spoke $s$ (inclusive) and spoke $s+1$ (modulo $q$, exclusive).
This arrangement can be given a coloring in the following way:
color an edge in sector $s$ of $\{3,q\}'$ with color $k+s$ modulo $q$, where $k$ is the color of the 
corresponding edge in sector $0$ of $\{4,q-2\}$.

To show that the coloring just described is a shortest-path-invariant partial addressing scheme, it will suffice to verify the hypotheses of Lemma 2.8.
These hypotheses constrain the coloring of edges incident with a common vertex, or incident with a common face.
It is clear that the edges incident with the origin receive distinct colors, so it will suffice to consider vertices that lie in a sector.
If all the edges incident with a vertex or face lie within one sector, then the hypothesis holds
automatically, since the coloring of a sector in $\{3,q\}'$ is the same as that of sector $0$ in 
$\{4,q-2\}$, except for a one-to-one renaming of the colors.
Thus it remains to consider vertices and faces that lie near spokes, in the sense that their incident edges lie in two consecutive sectors.
It will suffice to consider vertices and faces that lie near spoke $0$, since what happens near other spokes differs only by a one-to-one renaming of the colors.
For a vertex  to have some incident edges in sector $0$ and some in sector $q-1$, it must lie on spoke $0$.
It coloring in $\{3,q\}'$ differs from its coloring in $\{4,q-2\}$ only in that some edges that are colored $q-3$ in $\{4,q-2\}$ are colored $q-1$ in $\{3,q\}'$.
Since the color $q-1$ is not used in $\{4,q-2\}$, this cannot create a conflict.

To verify the hypotheses for faces, we consider two cases according to the parity of $q$.
If $q = 2r$ is even, then along spoke $0$ there is an edge joining each consecutive pair of vertices, and these edges are alternately colored $0$ and $r-1$.
The edge on spoke $0$ colored $0$, have their consecutive edges of their incident faces in sector $q-1$ colored $q-1$, and their opposite edge is colored $0$, satisfying the hypothesis.
Similarly, the edges on spoke $0$ colored $r-1$, have their consecutive edges of their incident faces in sector $q-1$ colored $r$, and their opposite edge  is colored $r-1$, satisfying the hypothesis.

If $q=2r+1$ is odd, then along spoke $0$ there is a periodic pattern with period six:
edges along the spoke alternate with faces that straddle the spoke, giving a graph-theoretic period of three (since an edge spans one generation, while a face spans two). 
But since ``clockwise'' and ``counterclockwise'' alternate generations, the pattern of the coloring repeats only after six generations.
The edges that lie along spoke $0$ are all colored $0$, their consecutive edges in their incident faces in sector $q-1$ are colored $q-1$ or $1$, and their opposite edges are colored $0$, satisfying the hypothesis.
The faces that straddle spoke $0$ each have one pair of opposite edges colored $r$ and another pair colored $r+1$, again satisfying the hypothesis.
\QED
\sk

\heading{3. Results for Transient Faults}

In this section, we shall establish the results shown in table (1.1) pertaining to transient faults.
Specifically, we shall show that a tessellation $\{p,q\}$ cannot tolerate transient faults if $q=2$,
if $q\le 4$ and $p<\infty$, or if $q\le 6$ and $p=3$.
The tessellation $\{p,q\}$ can in fact tolerate transient faults in all other cases, but in this section,
we shall only prove this fact for $p=\infty$ and $q\ge 3$,
for $p=4$ and $q\ge 5$,
and for $p=3$ and $q\ge 7$, because in all the remaining cases, $5\le p < \infty$,
the tessellation can even tolerate combined faults, as will be shown in the following section.

All of our ultimate results concern binary cellular automata that use symmetric majority voting at each time step.
These automata possess many useful properties.
Firstly, they are {\it time-invariant}; that is, the transition function is unchanged from time step to time step.
Secondly, they are {\it homogeneous\/} and {\it isotropic}; that is, the transition function is the same at each cell, and is unchanged by
rotation around a vertex or face, or by reflection across an edge (exchanging its two incident faces) or along an edge (exchanging its two incident vertices).
This property follows from the fact that the automorphism group of the underlying tessellation acts transitively on the vertices and oriented edges of its graph (see Coxeter and Moser [C2]).
Finally, the transition function is {\it monotone\/} (changing an input from $0$ to $1$ cannot change the output from $1$ to $0$) and
{\it self-dual\/} (complementing all the inputs complements the output).
In this section we will also need to consider automata (``sped-up'' automata) that do not use symmetric majority voting (but rather use an iterated version of it); these automata will, however, possess all the properties just enumerated.
But in the following section, we will need to consider automata (``weakened'' automata) that are not homogeneous or isotropic, and whose transition functions are not self-dual.

Our main tool for positive results concerning transient faults is Toom's theorem (Theorem 3.2 below).
As stated by Toom [T2], this theorem does not require time-invariance, montonicity, or self-duality.
Our statement of it below, however, has been simplified by assuming these properties,
which hold in all the cases we will consider.
It has also been simplified by replacing Toom's definition of transient fault tolerance, which involves taking a supremum over a set of probability measures, by the definition given above in Section 1, which involves an adversary whose possible actions correspond to the probability measures in that set.

Let $A$ be a monotone and self-dual binary automaton.
Since $A$ is self-dual, to determine whether it is fault-tolerant it will suffice to consider the case in which it is trying to remember the bit $x=0$, so that all cells are initially in state $0$, and a cell is in error at a particular time if it is in state $1$ at that time.
The level of fault tolerance in this case will then also hold in the dual case $x=1$.
Since $A$ is monotone, the adversary has a simple optimal strategy:
take any opportunity provided by a fault to put a cell into state $1$ at a particular time.
(This strategy is optimal because following it cannot foreclose any opportunity to put any cell 
into state $1$ at any other time.)

We shall now introduce a graph, the ``dependence graph'', that embodies information about the transition function of a cellular automaton (as well as about its underlying tessellation).
Let $a$ be a cell in $A$.
We shall say that a set $S$ of cells in $A$ is an {\it error set\/} for $a$ if the transition function for 
$a$ assumes the value $1$ whenever the states of all the cells in $S$ assume the value $1$.
We shall say that an error set $S$ for $a$ is {\it minimal\/} if no proper subset of $S$ is an error set for $a$.
(For a monotone transition function, the minimal error sets uniquely determine the transition function, since they correspond to the ``min-terms'' of that function.)
We shall say that cell $a$ {\it depends\/} on cell $b$ (or that $a$ is a {\it dependent\/} of $b$)
if $b$ belongs to some minimal error set for $a$.
We shall define the {\it dependence graph\/} $\Ga = (V, H)$ of $A$ to be the directed graph that has the same set of vertices $V$ as the graph $G = (V,E)$ of the underlying tessellation, but that has a directed edge $(a,b)\in H$ when and only when $a$ is a dependent of $b$.
If $A$ is the automaton using symmetric majority voting in the tessellation $\{p,q\}$,
then $(a,b)\in H$ if $\{a,b\}\in E$, or if $a=b$ and $q$ is even.
Every vertex in $\Ga$ has both in-degree and out-degree both equal to either $q$ or $q+1$, whichever is odd.

We shall say that a cell $a$ is a {\it guardian\/} of cell $b$ if $b$ is a dependent of $a$.
All of the automata with which we shall deal will be ``reciprocal'', in the sense that $a$ is a dependent of $b$ if and only if $b$ is a dependent of $a$.
(The terms ``dependent'' and ``guardian'' will thus be equivalent, and we could if we wished use them interchangeably, or use one to the exclusion of the other.
We shall, however, say that ``$a$ is a dependent of $b$'' to emphasize that the state of $a$ at time $t+1$ depends on that of $b$ at time $t$, and that ``$a$ is a guardian of $b$'' to emphasize that the state of $b$ at time $t+1$ depends on that of $a$ at time $t$.)

The negative results of this section use the well known notion of a ``self-sustaining island'' of errors.
A {\it self-sustaining island\/} in a tessellation is a finite set $I$ of vertices with the following property:
if the cells at vertices in $I$ are in error at some time $t_0$, while all other cells are not in error, then the cells in $I$ will remain in error for all times $t\ge t_0$ (because each cell in $I$ has a majority
of its guardians also in $I$).
Clearly, if a homogeneous and monotone cellular automaton has a self-sustaining island, then it cannot tolerate transient errors.
For by homogeneity, every cell is contained in a translated version of that island.
Sooner or later, every cell of that version will suffer a fault simultaneously, and from that time on,
every cell of the version will remain in error.
Thus, every cell of the tessellation is eventually in error, and remains in error for all subsequent time.

If $q=2$, we may take any two adjacent vertices as a self-sustaining island: each vertex of the island has a neighbor in the island and this, together with the vertex itself, gives the required majority of two out of three votes to keep the island in error.

If $3\le q\le 4$ and $p<\infty$, we may take the $p$ vertices around a face as a self-sustaining island:
each vertex has two neighbors in the island and thee neighbors, together with the vertex itself if $q=4$,
give the required two out of three, or three out of five if $q=4$, votes required to keep the island in error.

Finally, if $5\le q\le 6$ and $p=3$, we may take a vertex together with its $q$ neighbors as a 
self-sustaining island: each vertex now has three at least neighbors in this island, and these neighbors, together with the vertex itself if $q=6$, give the required three out of five, or four out of seven if $q=6$, votes required to keep the island in error.
These arguments yield the following theorem.

\label{Theorem 3.1:}
The tessellation $\{p,q\}$ cannot tolerate transient faults if $q=2$,
if $q\le 4$ and $p<\infty$, or if $q\le 6$ and $p=3$.

We shall now introduce a graph, the ``extended dependence graph'',
that embodies temporal as well as spatial relationships among the states of various cells at various times.
We shall define the {\it extended dependence graph\/} of $A$ to be the directed graph
$\Ga^* = (V^*, H^*)$ that has as its set of {\it nodes\/} $V^* = V\times \bfN$ (where $\bfN = 
\{0, 1, 2, \ldots\}$ denotes the set of natural numbers), and directed {\it arcs\/}
$\((a,s), (b,t)\)\in H^*$ when and only when $(a,b)\in H$ and $s=t+1$
(that is, when the state of cell $a$ at time $s$ depends directly on that of cell $b$ at time $t$).
We shall say that a set $S^*$ of nodes in $\Ga^*$ is an {\it error set\/} (respectively, {\it minimal error set}) for the node $(a,t+1)$ if $S^* = S\times \{t\}$, where the set of vertices $S$ in $\Ga$ is an error set (respectively, minimal error set) for $a$.

Our positive results on transient fault tolerance will all use the following theorem.

\label{Theorem 3.2:}
(A.~L. Toom, 1980, [T2])
Let $\Ga^* = (V^*,H^*)$ be the  extended dependence graph of a time-invariant, monotone and self-dual binary cellular automaton $A$.
Suppose that for some $n$ there exist $n$ functions $L_1, L_2, \ldots, L_n : V^* \to \bfR$
( where $\bfR$ denotes the set of real numbers)
that satisfy the following three conditions.
\medskip

\item{(1)} 
There exists a constant $M$ such that, for every arc $(a,b) \in H^*$ and every $k$ ($1\le k\le n$), 
$\abs{L_k(b) - L_k(a)} \le M$.

\item{(2)}
For every $a\in V^*$, $\sum_{1\le k\le n} L_k(a) = 0$.

\item{(3)}
For every $a \in V^*$, every error set $S^*\subseteq V^*$ for $a$ and every $k$ 
($1\le k\le n$), there exists  $b\in S^*$ such that $L_k(b) - L_k(a) \ge 1$.
\medskip
\noindent Then $A$ tolerates transient faults.
Indeed, if the fault rate is $\ep>0$, the probability that a particular cell is in error at a particular time is $O(\ep)$.

\label{Lemma 3.3:}
Let $A$ be a cellular automaton with a shortest-path-invariant addressing scheme using $l$ colors.
Suppose that for every vertex $a$  in $A$ and for every error-set $S$ for $a$,
(i) there is a vertex $b\in S$ with $\norm{b} \ge \norm{a} + 1$, and 
(ii) for each color $k\in\{0,\ldots,l-1\}$, there is a vertex $b_k \in S$ with 
$\norm{b_k}_k \le \norm{a}_k$.
Then $A$ tolerates transient faults.

\label{Proof:}
We shall apply Theorem 3.2 to $\Ga^*$ with $n = l+1$.
Given a node $(a,t)$ in $\Ga^*$ and $k\in\{0,\ldots,l-1\}$, we define
$$L_{k+1}\((a,t)\) = -(l+1)\norm{a}_k - t,$$
and
$$L_{l+1}\((a,t)\) = (l+1)\norm{a} + l t.$$
Since $\norm{\cdots}$ and each $\norm{\cdots}_k$ change by at most $1$ when proceeding from 
a vertex in $G$ to one of its neighbors, it is clear that Condition (1) of Theorem 3.2 is satisfied with $M = l+2$.
Condition (2) is verified by the identity $\norm{a} = \norm{a}_0 + \cdots + \norm{a}_{l-1}$.
It remains to verify Condition (3).

Given vertex $a$, an error-set $S$ for $a$, and a color $k\in\{0,\ldots,l-1\}$, we can find
by hypothesis (ii) a vertex $b_k\in S$ such that $\norm{b_k}_k \le \norm{a}_k$.
We then have
$$\eqalign{
L_{k+1}\((b_k,t)\)
& = -(l+1)\norm{b_k} - t \cr
&\ge -(l+1)\norm{a} - t \cr
&= L_{k+1}\((a,t+1)\) + 1, \cr
}$$
which verifies Condition 3 for $L_1, \ldots, L_{n-1}$.
For $L_n$, we have by hypothesis (i) a vertex $b\in S$ such that $\norm{b} \ge \norm{a} + 1$.
We then have
$$\eqalign{
L_{l+1}\((b,t)\)
&= (l+1)\norm{b} + l t \cr
&\ge (l+1)\norm{a} + (l+1) + l t \cr
&= (l+1)\norm{a} + l (t+1) + 1 \cr
&= L_{l+1}\((a,t+1)\) + 1, \cr
}$$
which completes the verification of Condition 3.
\QED

For the positive results in this section, we shall begin with the simplest case, that of trees.

\label{Proposition 3.4:}
For $q\ge 3$ the tessellation $\{\infty,q\}$ tolerates transient faults.

\label{Proof:}
By Lemma 2.7, the tessellation $\{\infty, q\}$ has a shortest-path-invariant addressing scheme using $q$ colors, which we shall use to apply Lemma 3.3.

Consider a vertex $a$ in the graph of $\{\infty, q\}$ and an error-set $S$ for $a$.
We claim that $S$ must contain at least two neighbors of $a$ (that is, guardians of $a$
distinct from $a$ itself).
If $q = 2r+1$ is odd, then $S$ must contain at least $r+1$ guardians (which, since $q$ is odd, are also neighbors) of $a$,
and since $q\ge 3$, we have $r+1\ge 2$, verifying the claim for odd $q\ge 3$.
If $q = 2r$ is even, then $S$ must again contain contain $r+1$ guardians of $a$
(all but at most one of which are neighbors), and since $q\ge 4$, we have $r\ge 2$, verifying the claim for 
even $q\ge 4$.

We may now verify the hypotheses of Lemma 3.3.
We have at least two neighbors of $a$ in $S$.
At most one of these neighbors can be the parent of $a$, so at least one such neighbor $b$
must be a child of $a$ and thus
must satisfy $b\in S$ with $\norm{b} \ge \norm{a} + 1$.
For $k\in\{0,\ldots,l-1\}$, only one of the neighbors of $a$ can be joined to $a$ by an edge colored $k$,
so at least one neighbor $b_k$ must satisfy $\norm{b_k}_k \le \norm{a}_k$.
Thus Lemma 3.3 applies, completing the proof of the proposition.
\QED

Next we turn to the case $p=4$, with $q\ge 5$.

\label{Proposition 3.5:}
For $q\ge 5$ the tessellation $\{4,q\}$ tolerates transient faults.

\label{Proof:}
By Lemma 2.9, the tessellation $\{4,q\}$ has a shortest-path-invariant addressing scheme using $q$ colors, which we shall use to
apply Lemma 3.3.

Consider a vertex $a$ in the graph of $\{4,q\}$ and an error-set $S$ for $a$.
We claim that $S$ must contain at least three neighbors of $a$.
If $q = 2r+1$ is odd, then $S$ must contain at least $r+1$ guardians (which, since $q$ is odd, are also neighbors) of $a$,
and since $q\ge 5$, we have $r+1\ge 3$, verifying the claim for odd $q\ge 5$.
If $q = 2r$ is even, then $S$ must again contain contain $r+1$ guardians of $a$
(all but at most one of which are neighbors), and since $q\ge 6$, we have $r\ge 3$, verifying the claim for 
even $q\ge 4$.

We may now verify the hypotheses of Lemma 3.3.
We have at least three neighbors of $a$ in $S$.
At most two of these neighbors can be parents of $a$, so at least one such neighbor $b$
must be a child of $a$ and thus
must satisfy $b\in S$ with $\norm{b} \ge \norm{a} + 1$.
For $k\in\{0,\ldots,l-1\}$, only one of the neighbors of $a$ can be joined to $a$ by an edge colored $k$,
so at least one neighbor $b_k$ must satisfy $\norm{b_k}_k \le \norm{a}_k$.
Thus Lemma 3.3 applies, completing the proof of the proposition.
\QED

For the cases with $p=3$, the application of Toom's theorem presents a new difficulty.
In the applications for $p=\infty$ and $p=4$, an error at the origin was eventually corrected because each cell eventually received correct information from its children.
Consider, however, the automaton with symmetric majority voting on $\{3,7\}$.
For the cells with two parents as well as two siblings, these four neighbors constitute a majority of the seven votes that determine the next state.
Thus such cells may remain in an erroneous state in spite of all their children being in the correct state.
This situation is rescued by the fact that, according to Lemma 2.2, two such cells cannot be siblings of each other; that is, of any two siblings, at least one must be a single-parent cell that cannot be led into error unless at least one of its children is in error.
This circumstance, however, means that it may take two time steps for an error in a two-parent cell to be corrected by information from a ``niece'' or ``nephew'' (a child of a sibling).

Our strategy for $p=3$ will thus be as follows.
Let $A$ be the symmetric majority voting automaton on the tessellation $\{3,q\}$ with $q\ge 7$.
We shall consider a ``sped-up'' version $A^2$ of the symmetric majority voting automaton $A$:
the transition function of $A^2$ will accomplish in one time step what the transition function of $A$ accomplishes in two time steps.
Transferring attention to $A^2$ is equivalent to assuming that faults occur only at every second time step.
We shall then prove Lemma 3.6, a generalization of Lemma 3.3 that will allow us to deduce from a shortest-path-invariant partial addressing scheme for $A$ the conclusion that $A^2$ tolerates transient faults.
Lemma 2.10 provides a shortest-path-invariant partial addressing scheme for $A$, allowing us to conclude in Lemma 3.7
that $A^2$ tolerates transient faults.
In Lemma 3.8 we shall show that if $A^2$ tolerates transient faults, then so does $A$, allowing us to conclude in Proposition 3.9 that $A$ tolerates transient faults.

The strategy just described applies to situations in which the speed-up is by a factor greater than two, and in which the transition function of $A$ is not necessarily symmetric majority voting.
We shall therefore prove Lemmas 3.6 and 3.8 in greater generality than we need in this paper,
in the hope that they will find further applications.

If $A$ is a time-invariant cellular automaton, its {\it $\ka$-fold speed-up\/} $A^\ka$ (where $\ka\ge 2$ is an integer) is the the cellular automaton (with the same cells and the same set of states for each cell) whose state transition function does in one time step what that of $A$ does in $\ka$ time steps.
If $A$ is monotone and self-dual, then so is $A^\ka$.

\label{Lemma 3.6:}
Let $A$ be a time-invariant,  monotone and self-dual binary cellular automaton with a 
shortest-path-invariant partial addressing scheme using $l$ colors.
Let $B = A^\ka$ be the $\ka$-fold speed-up of $A$.
Suppose that for every vertex $a$  in $B$ and for every error-set $S$ for $a$ in $B$,
(i) there is a vertex $b\in S$ with $\norm{b} \ge \norm{a} + 1$, and 
(ii) for each color $k\in\{0,\ldots,l-1\}$, there is a vertex $b_k \in S$ with 
$\norm{b_k}_k \le \norm{a}_k$.
(Here $\norm{\cdots}$ and $\norm{\cdots}_k$ refer to path lengths in $A$ and the 
shortest-path-invariant partial addressing scheme of $A$.)
Then $B$ tolerates transient faults.

\label{Proof:}
Let $\De = (V, K)$ be the dependence graph of $B$, and let
$\De^* = (V^*, K^*)$ be the extended dependence graph of $B$.
We shall apply Theorem 3.2 to $\De^*$ with $n = l+1$.
Given a node $(a,t)$ in $\De^*$ and $k\in\{0,\ldots,l-1\}$, we define
$$L_{k+1}\((a,t)\) = -(l+1)\norm{a}_k -  t,$$
and
$$L_{l+1}\((a,t)\) = (l+1)\norm{a} + l  t.$$
Since $\norm{\cdots}$ and each $\norm{\cdots}_k$ change by at most $\ka$ when proceeding from 
a vertex in $\De$ to one of its neighbors, it is clear that Condition (1) of Theorem 3.2 is satisfied with $M = \ka(l+2)$.
Condition (2) is verified by the identity $\norm{a} = \norm{a}_0 + \cdots + \norm{a}_{l-1}$.
It remains to verify Condition 3.

Given vertex $a$, an error-set $S$ in $B$ for $a$, and a color $k\in\{0,\ldots,l-1\}$, we can find
by hypothesis (ii) a vertex $b_k\in S$ such that $\norm{b_k}_k \le \norm{a}_k$.
We then have
$$\eqalign{
L_{k+1}\((b_k,t)\)
& = -(l+1)\norm{b_k} -  t \cr
&\ge -(l+1)\norm{a} -  t \cr
&= -(l+1)\norm{a} - (t+1) + 1 \cr
&= L_{k+1}\((a,t+1)\) + 1, \cr
}$$
which verifies Condition 3 for $L_1, \ldots, L_{n-1}$.
For $L_n$, we have by hypothesis (i) a vertex $b\in S$ such that $\norm{b} \ge \norm{a} + 1$.
We then have
$$\eqalign{
L_{l+1}\((b,t)\)
&= (l+1)\norm{b} + l  t \cr
&\ge (l+1)\norm{a} + (l+1) + l t \cr
&= (l+1)\norm{a} + l (t+1) + 1 \cr
&= L_{l+1}\((a,t+1)\) + 1, \cr
}$$
which completes the verification of Condition 3.
\QED

\label{Lemma 3.7:}
Let $A$ be the binary automaton using symmetric majority voting on $\{3,q\}$ with $q\ge 7$,
and let $B = A^2$.
Then $B$ tolerates transient faults.

\label{Proof:}
By Lemma 2.10, $A$ has a shortest-path-invariant partial addressing scheme using $q$ colors.
We shall apply Lemma 3.6 with $l=q$ and $\ka=2$.

Consider a vertex $a$ and an error-set $S$ for $a$ in $B$.
There must exist an error-set $T$ for $a$ in $A$, and for each vertex $c$ in $T$, $S$ must include an error-set $S_c$ for $c$ in $A$.
We claim that $T$ must contain at least four guardians (which, since $q$ is odd, are neighbors) of $a$.
If $q = 2r+1$ is odd, then $T$ must contain at least $r+1$ neighbors of $a$,
and since $q\ge 7$, we have $s+1\ge 4$, verifying the claim for odd $q\ge 7$.
If $q = 2r$ is even, then $T$ must again contain contain $r+1$ guardians of $a$
(all but at most one of which are neighbors), and since $q\ge 8$, we have $s\ge 4$, verifying the claim for 
even $q\ge 6$.
By a similar argument, for each $c\in T$, $S_c$ must contain at least four neighbors of $c$.

We may now verify the hypotheses of Lemma 3.6.
To satisfy condition (i), we must take a step from $a$ to a neighbor $b\in S$ of $a$ in $B$ such that
$\norm{b} \ge \norm{a} + 1$.
To do this we shall take two steps in $A$, the first from $a$ to a neighbor $c\in T$ of $a$ in $A$, and the second from $c$ to a neighbor $b\in S_c$ of $c$ in $A$.
At each step we shall move from a vertex to a child of that vertex if possible, and to a sibling of that vertex if no move to a child is possible; it will never be necessary to move to a parent.
If a vertex has just one parent, then at most three of the neighbors in an error-set for that vertex are  parents or siblings, and thus a move to a child must be possible.
Only if a vertex has two parents can all four of the neighbors in an error-set for that vertex be parents or siblings, forcing a move to a sibling.

By Lemma 2.2,  two consecutive vertices in a generation cannot both be two-parent vertices.
Thus, of the two steps to neighbors in $A$,  neither can be to a parent, and at least one must be to a child, since two consecutive steps to siblings could only occur from two consecutive 
two-parent vertices in a generation, which is impossible.
It follows that neither step decreases the norm $\norm{\cdots}$, and at least one of the steps increases it by  $1$.
Thus these steps arrive at a vertex $b\in S$ with $\norm{b} \ge \norm{a} + 1$.

To satisfy condition (ii) for a given $k\in\{0,\ldots,l-1\}$, we must take a step from $a$ to a neighbor 
$b\in S$ of $a$ in $B$ such that $\norm{b}_k \le \norm{a}_k$.
To do this we shall again take two steps in $A$, the first from $a$ to a neighbor $c\in T$ of $a$ in $A$, and the second from $c$ to a neighbor $b\in S_c$ of $c$ in $A$.
We shall show that it is always possible to take such a step in a {\it safe\/} way, that is, without
increasing the $k$-th norm $\norm{\cdots}_k$.
Of the four or more neighbors in an error-set for a vertex, at most one can be a child
with greater $k$-th norm (a child connected by an edge with color $k$), and at most two can be siblings.
Thus at least one of the neighbors in an error-set must be either a parent or a child with
no greater $k$-th norm, to which a safe move is possible.
This completes the verification of hypothesis (ii) of Lemma 3.6, and thus completes the proof of the proposition.
\QED

\label{Lemma 3.8:}
Let $A$ be a time-invariant, homogeneous and monotone binary cellular automaton, and let 
$B = A^\ka$ be its $\ka$-fold speed-up.
Then if $B$ tolerates transient faults, so does $A$.

\label{Proof:}
We shall construct a coupling between the probability space of $A$ and that of $B$.
That is, we shall define on a common probability space a copy of the probability space of $A$
and a copy of the probability space of $B$.
This construction will be such that whenever the faults in $A$ are such as to allow the adversary 
to force a cell $v$ at a time $\ka t$ ($t\ge 1$) into error, the faults in $B$ will allow the adversary
to force the cell $v$ at time $t$ into error.
It will follow by induction on $t$  that the errors in $A$ at time $\ka t$ are stochastically dominated by those in $B$ at time $t$.
Since the probability of error at a cell in a monotone binary cellular automaton is a non-decreasing function of time, it will follow that the supremum of the probability of error in $A$ over time is bounded by that of $B$, so that if $B$ tolerates transient faults, then $A$ does as well.

Let $\Ga=(V,H)$ denote the dependence graph of $A$.
Let $\Ga_\ka = (V,H_\ka)$ denote the directed graph that has the same set of vertices $V$ as $\Ga$, but that has a directed edge $(v,w)\in H_\ka$ when and only when there is a path of length at most $\ka$ from $v$ to $w$ in $\Ga$.
Since $A$ is homogeneous, every vertex has the same in-degree and out-degree in $\Ga_\ka$.
Let $\la$ denote this common in- and out-degree.

Let $X_{(v,w),t}$ for $(v,w)\in H_\ka$ and $t\ge 0$ be independent Bernoulli random variables
with $\Pr[X_{(v,w),t} = 1] = \xi$.
Define $Y_{v,t}$ for $v\in V$ and $t\ge 0$ by
$$Y_{v,t} = \bigwedge_{(v,w)\in H_\ka} X_{(v,w),t}.$$
Since the $Y_{v,t}$ are defined from disjoint sets of independent random variables, they  are independent Bernoulli random variables with $\Pr[Y_{v,t} = 1] = \et$,
where $\et = \xi^\la$.
We may thus take $Y_{v,t}$ to be the indicator of a fault at cell $v$ at time $t$ in automaton $A$,
with fault rate $\et$.
Define $Z_{w,s}$ for $w\in V$ and $s\ge 0$ by
$$Z_{w,s} = \bigvee_{(v,w)\in H_\ka \atop 0\le \mu\le \ka-1} X_{(v,w),\ka s + \mu}.$$
Since the $Z_{w,s}$ are defined from disjoint sets of independent random variables, they  are independent Bernoulli random variables with $\Pr[Z_{w,s} = 1] = \ze$,
where $\ze = 1 - (1 - \xi)^{\ka\la}$.
We may thus take $Z_{w,s}$ to be the indicator of a fault at cell $v$ at time $\ka t$ in automaton 
$B$, with fault rate $\ze$.

Now if in automaton $A$ no faults occur at cell $v$ or any cell with in distance $\ka$ of $v$
at any time in the interval $[\ka t, \ka(t-1)]$, then the state of cell $v$ at time $\ka t$ will be the same function
the states of the cells within distance $\ka$ of $v$ at time $\ka(t-1)$ as the transition function for automaton $B$.
But if  in automaton $A$ any fault occurs at cell $v$ or any cell with in distance $k$ of $v$
at any time in the interval $[\ka t, \ka(t-1)]$, then a fault will occur at cell $v$ at time $t$  in automaton $B$.
It follows that the errors in automaton $A$ at time $\ka t$ are stochastically dominated by those in
automaton $B$ at time $t$.
Thus if $\ze$ is sufficiently small that $B$ tolerates transient faults, the corresponding value of
$\et = \(1 - (1 - \ze)^{1/\ka\la}\)^\la$ will be sufficiently small that $A$ tolerates transient faults.
\QED

Lemma 3.7 and Lemma 3.8 together prove the following proposition.

\label{Proposition 3.9:}
For $q\ge 7$ the tessellation $\{3,q\}$ tolerates transient faults.

Propositions 3.4, 3.5 and 3.9 together prove the following theorem.

\label{Theorem 3.10:}
The tessellation $\{p,q\}$ tolerates transient faults if $p=\infty$ and $q\ge 3$,
if $p=4$ and $q\ge 5$, or if $p=3$ and $q\ge 7$.
\sk

\heading{4. Results for Combined Faults}

In this section, we shall establish the results shown in Table (1.1) pertaining to combined faults.
Specifically, we shall show that a tessellation $\{p,q\}$ can tolerate combined faults if $p\ge 5$ and $q\ge 5$,
if $p=4$ and $q\ge 7$, or if $p=3$ and $q\ge 9$.
The tessellation $\{p,q\}$ cannot in fact tolerate combined faults in any other cases, but in this section
we shall only prove this for  $p=\infty$ and $3\le q\le 4$, for $p=4$ and $5\le q\le 6$, and for $p=3$ and 
$7\le q \le 8$, because in all the remaining cases the tessellation cannot even tolerate transient faults, as was shown in the preceding section.

For the negative results of this section, we shall introduce the notion of  a ``pier-supported bridge''.
A {\it pier-supported bridge\/}  in a tessellation is a pair $(I,J)$ of finite sets of vertices, where $I$ is called the {\it bridge\/} and $J$ is called the set of {\it piers}, with the following property: if the cells at the vertices
in $J$ are in error for all time, and the cells at the vertices in $I$ are in error at some time $t_0$,
then the cells in $I$ will remain in error for all time $t\ge t_0$ (because each cell in $I$ has at least
one-half of its guardians in $I\cup J$).
If the piers of a pier-supported bridge all suffer permanent faults, then eventually all the cells of the bridge will simultaneously suffer transient faults, and from that time onward all cells in the bridge
and its piers will be in error.
Let $v$ be a cell in a tessellation $G$.
If with probability $1$ there exists a pier-supported bridge whose piers all suffer permanent failures and whose bridge contains $v$, then $G$ does not tolerate combined faults.

Suppose first that $p=\infty$ and  $3\le q\le 4$.
Let $v_0$ be a cell, and let 
$\ldots, v_{-1}, v_0, v_1, \ldots$ be an indexing of the vertices around a face including 
$v_0$ such that $v_k$ and $v_{k+1}$ are joined by an edge for all $k\in\bfZ$.
Then for any $m\ge 1$ and $n\ge 1$, $(I_{m,n},J_{m,n})$, where $I_{m,n}=\{v_{-(m-1)}, \ldots, v_{n-1}\}$ and 
$J_{m,n}=\{v_{-m}, v_n\}$, form a pier-supported bridge, since each vertex in $I$ has two neighbors
(plus itself, if $q=4$) in $I_{m,n}\cup J_{m,n}$.
With probability $1$, there exists an $m\ge 1$ such that $v_{-m}$ suffers a permanent failure,
since this event is the union of countably many independent events, each of which occurs
with the same strictly positive probability.
Similarly, with probability $1$ there exists an $n\ge 1$ such that $v_n$ suffers a permanent failure.
Since these two events each occur with probability $1$, their intersection also occurs with probability 
$1$.
For this $m$ and $n$,  $(I_{m,n},J_{m,n})$ then forms a pier-supported bridge, and 
thus tessellations with $3\le q\le 4$ do not tolerate combined faults.

For the following we shall need the notion of a an ``\oes''.
If $e$ is an edge in a tessellation in which every face has degree four, the {\it \oes\/} containing 
$e$ is the smallest set $F$ of edges such that (1) $e\in F$ and (2) if $f\in F$ and $g$ is the edge opposite $e$ on either of the faces bounded by $f$, then $g\in F$.

Suppose next that $p=4$ and $5\le q\le 6$.
Let $v_0$ be a cell, let $e_0$ be an edge containing $v_0$, and let 
$F$ be the \oes{} containing $e_0$.
If $F$ were finite, the vertices bounding the edges contained in $F$ would constitute
a self-sustaining island, and the tessellation would not even tolerate transient errors.
(This case cannot arise, since it would contradict Proposition 3.5.)
Thus we may assume that $F$ is infinite.
Let $\ldots, e_{-1}, e_0, e_1, \ldots$ be an indexing of the edges in $F$ such that $e_k$ and $e_{k+1}$
are opposite edges of a face for all $k\in \bfZ$.
Then for any $m\ge 1$ and $n\ge 1$, $(I_{m,n},J_{m,n})$, where $I_{m,n}=e_{-(m-1)}\cup\cdots\cup e_{n-1}$ and 
$J_{m,n}=e_{-m}\cup e_n$, form a pier-supported bridge, since each vertex in $I$ has three neighbors
(plus itself, if $q=6$) in $I_{m,n}\cup J_{m,n}$.
With probability $1$, there exists an $m\ge 1$ such that both vertices of $e_{-m}$ suffer permanent failures,
since this event is the union of countably many independent events, each of which occurs
with the same strictly positive probability.
Similarly, with probability $1$ there exists an $n\ge 1$ such that both vertices of $e_n$ suffer permanent failures.
Since these two events each occur with probability $1$, their intersection also occurs with probability 
$1$.
For this $m$ and $n$,  $(I_{m,n},J_{m,n})$ then forms a pier-supported bridge, and 
thus tessellations with $p=4$ and $5\le q\le 6$ do not tolerate combined faults.

Suppose finally that $p=3$ and $7\le q\le 8$.
Let $v_0$ be a vertex in $\{3,q\}$
We consider the tessellation $\{3,q\}'$, obtained from $\{3,q\}$ by deleting all sibling edges with respect to the origin $v_0$.
Every face of $\{3,q\}'$ has degree four, so the argument of the previous paragraph shows that with probability $1$, there is a pier-supported bridge $(I,J)$ with $v_0\in I$ and permanent faults at every vertex in $J$.
Thus the tessellations with $p=3$ and $7\le q\le 8$ do not tolerate combined faults.
These arguments yield the following theorem.

\label{Theorem 4.1:}
The tessellation $\{p,q\}$ cannot tolerate combined faults if $p=\infty$ and $3\le q\le 4$, if 
$p=4$ and $5\le q\le 6$, or if  $p=3$ and $7\le q\le 8$.

All of our positive results in Section 2 were derived from Toom's theorem (Theorem 3.2).
That theorem applies to transient faults, but not to combined faults.
We shall next prove an analogous theorem that applies to combined faults.

Let $A$ be a binary cellular automaton using symmetric majority voting.
Because majority voting is self-dual, we may assume without loss of generality that $A$
is trying to remember the bit $0$, so that all cells are initially in state $0$.
A cell is in error if it is in state $1$, and because majority voting is monotone, the optimal strategy
for the adversary is to take any opportunity to put any cell in state $1$.

To show that a cellular automaton $A$ with majority voting tolerates combined faults, we shall
use a monotone binary automaton $B$ that will be called a {\it weakening\/} of $A$.
The automaton $B$ will be obtained from $A$ by specifying, for each cell $a$ in $A$, a set $I(a)$ of guardians of $a$ to be {\it ignored}.
The cells of $B$ will be the same as the cells of $A$, but the transition function of a cell $a$ in $B$
will be obtained from that of the cell $a$ in $A$ by substituting $1$ for each argument corresponding to an ignored cell in $I(a)$.
Because of the monotonicity of the transition functions, for any disposition of permanent and transient faults in space and time, the errors in $B$ will stochastically dominate the errors at corresponding places and times in $A$.
Thus the probability that cell $a$ of $A$ is in error at time $t$ is at most the probability that cell $a$ of $B$ is in error at time $t$.

The transition function for a cell in $A$ is a threshold function, which assumes the value $1$
if the number of $1$s among its arguments exceeds a certain number called the {\it threshold\/}
of the cell.
(For majority voting, the threshold of every cell is $r+1$ if the cell has $2r+1$ guardians, that is, if it has an even number $2r$ or an odd number $2r+1$ of neighbors.)
The transition function for a cell in $B$ is also a threshold function, whose threshold, called the 
{\it reduced threshold}, is obtained by reducing the original threshold by the number of vertices ignored
by the cell.

To bound the probability that cell $a$ of $B$ is in error at time $t$, we shall consider a 
particular disposition of permanent and transient faults in space and time for which cell $a$ is in error at time $t$, and construct a graph $\Xi^* = (X^* \cup Y^*, Z^*)$ called the {\it extended explanation graph\/} rooted at cell $a$ at time $t$, which will be a finite subgraph of the
extended dependence graph of $\De^* = (V^*, K^*)$ of $B$.
The set $X^* \cup Y^*$  of nodes of $\Xi^*$ is partitioned into two disjoint subsets,
$X^* \cap Y^* = \emptyset$: the nodes in $X^*$ will be called {\it terminal\/} nodes,
while those in $Y^*$ will be called {\it non-terminal\/} nodes.

We shall define the extended explanation graph $\Xi^*$ by giving an algorithm for constructing it.
This algorithm will work in a breadth-first fashion, and will maintain a queue of nodes, which will initially be empty.
If a node $(b,s)$ appears in the queue at any time during the execution of the algorithm,
the cell $b$ is in error at time $s$.
The sets $X^*$, $Y^*$ and $Z^*$ will be regarded as variables that are initially 
empty, that may have elements added to them from time to time, and whose final values will define the extended explanation graph $\Xi^*$.
\sk

\item{(1)}
Put the node $(a,t)$ into the queue.

\item{(2)}
While the queue is non-empty, perform the following sequence of steps until instructed to stop.

\itemitem{(2.1)}
Dequeue a node $(b,s)$ from the queue.

\itemitem{(2.2)}
If $X^* \cup Y^*$ currently contains a node $(b,r)$ for any $s<r\le t$, stop.

\itemitem{(2.3)}
If the cell $b$ is at fault at time $s$ (either because of a transient fault at $b$ at time $s$, or because of a permanent fault at $b$), put the node $(b,s)$ into $X^*$
(so that $(b,s)$ becomes a terminal node); then stop.

\itemitem{(2.4)}
Put the node $(b,s)$ into $Y^*$ (so that $(b,s)$ becomes a non-terminal node).
Since the cell $b$ is in error, but not at fault, at time $s$, there must exist an error set $S^*$
for $(b,s)$ in $\De^*$ such that (i) the number of nodes in $S^*$ is at least the reduced threshold of $b$, and (ii) for each node $(c,r)$ in $S^*$, $s=r+1$ and the cell $c$ is in error at time $r$.
For each node $(c,r)$ in $S^*$, put the arc $\((b,s),(c,r)\)$ into $Z^*$ and put the node $(c,r)$ into the queue; then stop.
\medskip
It is clear that a cell $b$ can appear in a node $(b,s)$ in $X^* \cup Y^*$ for at most one value of 
$s$, and that a node $(b,s)$ can appear in $X^*$ or $Y^*$ but not both.
The out-degree of a terminal node is zero, and the out-degree of a non-terminal node $(b,s)$ is a least the reduced threshold of $b$.

From the extended explanation graph $\Xi^*$ just constructed, we shall construct the 
{\it explanation graph\/} $\Xi = (X\cup Y,Z)$, which is a finite subgraph of the dependency graph 
$\De = (V, K)$ of $B$, by
``projecting out the time'':
The graph $\Xi$ contains a terminal vertex $b$ if $\Xi^*$ contains a terminal node $(b,s)$ for any time $s$, contains a non-terminal vertex $b$ if $\Xi^*$ contains a non-terminal node $(b,s)$ for any time $s$, and contains an edge $(b,c)$ if $X^*$ contains an edge $\((b,s),(c,r)\)$ for any times $s$ and $r$.
It is clear that a vertex can appear in $X$ or $Y$, but not both.
The out-degree of a terminal vertex is zero, and the out-degree of a non-terminal vertex  $b$ is a least the reduced threshold of $b$.

The construction of $\Xi^*$ identifies for each terminal vertex $b$ of $\Xi$ a particular time 
(the time $s$ for which $(b,s)$ is a node in $\Xi^*$) at which $b$ is at fault.
The probability that a particular cell is at fault at a particular time is at most
$\ep = 1 - (1-\al)(1-\be)$, and these probabilities are independent for distinct cells.
Thus the probability that a particular subgraph $\Xi$ of the dependency  graph $\De$ of $B$
is the explanation graph for a cell $a$ that is in error at time $t$ is at most $\ep^m$, where $m$ is the number of terminal vertices in $\Xi$.

An explanation graph rooted at $a$ and containing $n$ vertices is an acyclic directed graph with source $a$, and thus
can be described by the sequence of edge traversals in a depth-first search.
Each vertex has total degree (in-degree plus out-degree) at most $q$, so if there are $n$ vertices,
there are at most $2qn$ edge traversals in the sequence.
Thus there are at most 
$$\sum_{1\le k\le 2qn} q^k = {q^{2qn+1} - 1 \over q - 1} \le q^{2qn+1}$$
possible explanation graphs rooted at $a$ and having $n$ vertices.

\label{Theorem 4.2:}
Let  $A$ be a time-invariant monotone binary automaton.
Suppose that for each cell $a$ in $A$ there exists a 
weakening $B$ of $A$ such that there exists a constant $M$ such that, for every explanation graph $\Xi$ rooted at $a$, the total number of vertices in $\Xi$ is at most $M$ times the number of terminal vertices in $\Xi$.
Then $A$ tolerates combined faults.

\label{Proof:}
If cell $a$ is in error at some time $t$, faults must occur corresponding to the terminal vertices of an explanation graph rooted at $a$ at time $t$.
Since there are at most $q^{2qn+1}$ possible explanation graphs with $n$ vertices, there are at most $q^{2qMm+1}$ possible explanation graphs with $m$ terminals.
Thus we have
$$\eqalign{
\Pr[a\hbox{\ in error at time\ }t]
&\le q \sum_{m\ge 1} q^{2qMm} \ep^{m} \cr
 &\le {q^{2qM+1} \ep \over 1 - q^{2qM} \ep}. \cr
 }$$
 This bound tends to zero as $O(\ep)$ as $\ep$ tends to zero, which completes the proof.
 \QED

We shall apply Theorem 4.2 in a sequence of results of increasing intricacy.
Our applications of Theorem 4.2 will use a ``balance of payments'' argument of the following form.
We start by putting a certain positive amount of money on each non-terminal vertex of the 
explanation graph, but no money on terminal vertices.
We shall then move a certain amount of money across each edge, taking it from the source vertex of the edge and giving it to the target vertex.
These movements will be chosen so as to leave no money on any non-terminal vertex, so that
(since money is conserved in the movements) all the money is now on the terminal vertices.
If $r$ is a lower bound to the amount of money initially on any non-terminal vertex, and $s$ is an upper bound to the amount of money on any terminal vertex, then $r$ times the number of non-terminals is at most $s$ times the number of terminals.
Thus we may take $M = 1 + s/r$ in Theorem 4.2.

Our applications of the balance of payments argument will use an initial distribution of money 
determined in the following way.
We will associate with each edge in the explanation graph a certain amount of money to be 
moved across it.
We then put on each non-terminal vertex the sum of the amounts associated with the edges directed out of the vertex, minus the sum of the amounts associated with the edges directed into the vertex.
This rule ensures that after the movement of funds across edges, each non-terminal is left with no money.

\label{Proposition 4.3:}
For $q\ge 5$, the tessellation $\{\infty,q\}$ tolerates combined faults.

\label{Proof:}
For $q\ge 5$, let $A$ be the cellular automaton based on the tessellation $\{\infty,q\}$ with symmetric majority voting.
Let $a$ be a cell in $A$, and take $a$ as the origin of the tessellation.
We define the weakened automaton $B$ by having each cell $b$ ignore (1) its parent
(unless $b$ is the origin) and (2) itself (if $q$ is even).
Since $q$ is at least $5$ ($6$ if $q$ is even), the threshold of each vertex is at least $3$ ($4$ if $q$ is even).
Thus the reduced threshold of each vertex is at least $2$.

The explanation graph is a finite tree with $a$ as its root and with the terminal vertices as its leaves.
We associate a flow of $1$ dollar with each edge.
Every non-terminal has in-degree at most $1$ (from its parent), and out-degree at least $2$ (its reduced threshold), and thus initially has at least $1$ dollar.
Every terminal has in-degree at most $1$, and thus finally has at most $1$ dollar.
We conclude that the number of terminals is at least the number of non-terminals, and thus that we may take $M=1+1/1=2$ in Theorem 4.2.
\QED

\label{Proposition 4.4:}
For $q\ge 7$, the tessellation $\{4,q\}$ tolerates combined faults.

\label{Proof:}
For $q\ge 7$, let $A$ be the cellular automaton based on the tessellation $\{4,q\}$ with symmetric majority voting.
Let $a$ be a cell in $A$, and take $a$ as the origin of the tessellation.
We define the weakened automaton $B$ by having each cell $b$ ignore (1) its parent or parents
(unless $b$ is the origin) and (2) itself (if $q$ is even).
Since $q$ is at least $7$ ($8$ if $q$ is even), the threshold of each vertex is at least $4$ ($5$ if 
$q$ is even).
Thus the reduced threshold of each one-parent vertex is at least $3$, and the reduced threshold of each two-parent vertex is at least $2$.
If we were to associate a flow of $1$ dollar with each edge, two-parent non-terminal vertices could have both in-flow $2$ (from their parents) and out-flow $2$ (their reduced thresholds), and thus initially have $0$ dollars.
To ensure that every non-terminal initially has strictly positive funds, we shall distinguish two types of  edges in the explanation graph, with different flows for each type.

Say that an edge is {\it special\/} if it is directed from a one-parent parent to a two-parent child.
We associate a flow of $1$ dollar with each special edge, and a flow of $3$ dollars with each other edge.

A one-parent vertex has in-flow at most $3$ (from its parent), since an in-edge to a one-parent vertex cannot be special.
It has out-degree at least $3$ (its reduced threshold).
Of these out-edges, at most two can be special (since only a leftmost or rightmost out-edge can be special), so the outflow is at least $1+1+3 = 5$.
Thus every one-parent vertex initially has at least $5-3=2$ dollars.

A two-parent vertex has in-degree $2$ (from its parents), but at least one of these parents must be a one-parent vertex (two two-parent parents would be consecutive, with no intervening one-parent vertices, contradicting Lemma 2.1).
Thus the in-flow is at most $1+3 = 4$.
It has out-degree at least $2$ (its reduced threshold), and thus its out-flow is at least $3+3 = 6$
(since an out-edge from a two-parent vertex cannot be special).
Thus every two-parent vertex initially has at least $6-4=2$ dollars.

A terminal vertex can have in-degree at most $2$, thus in-flow at most $3+3=6$, and thus finally at most $6$ dollars.
We conclude that $6$ times the number of terminals is at least $2$ times the number of non-terminals, and thus that we may take $M=1+6/2=4$ in Theorem 4.2.
\QED

\label{Proposition 4.5:}
For $q\ge 9$, the tessellation $\{3,q\}$ tolerates combined faults.

\label{Proof:}
For $q\ge 9$, let $A$ be the cellular automaton based on the tessellation $\{3,q\}$ with symmetric majority voting.
Let $a$ be a cell in $A$, and take $a$ as the origin of the tessellation.
We define the weakened automaton $B$ by having each cell $b$ ignore (1) its parent or parents
(unless $b$ is the origin), (2) its siblings (if it is a one-parent vertex), and (3) itself (if $q$ is even).
Since $q$ is at least $9$ ($10$ if $q$ is even), the threshold of each vertex is at least $5$ ($6$ if $q$ is even).
Thus the reduced threshold of each one-parent vertex is at least $2$, and the reduced threshold of each two-parent vertex is at least $3$.

We associate a flow of $2$ dollars with each sibling edge, and a flow of $3$ dollars with each other edge. 
Of the two siblings of a one-parent vertex, at least one must be another one-parent vertex
(since otherwise there would be two two-parent vertices with only one intervening one-parent vertex, contradicting Lemma 2.2).
Thus a one-parent vertex can have at most one sibling in-edge.
Its total in-flow is thus $3$ (from its parent) plus at most $2$ (from a sibling), for a total of $5$.
Its out-degree is at least $2$ (its reduced threshold), so its out-flow is at least $3+3=6$.
Thus every one-parent vertex initially has at least $6-5=1$ dollar.

A two-parent vertex has in-degree $2$ (from its parents), and thus has in-flow $3+3=6$.
It has out-degree at least $3$ (its reduced threshold).
At most $2$ of these out-edges can be sibling edges, so its out-flow is at least $2+2+3=7$.
Thus every $2$-parent vertex initially has at least $7-6=1$ dollar.

A terminal vertex can have in-degree at most $2$, thus in-flow at most $3+3=6$, and thus finally at most $6$ dollars.
We conclude that $6$ times the number of terminals is at least  the number of non-terminals, and thus that we may take $M=1+6/1=7$ in Theorem 4.2.
\QED

\label{Proposition 4.6:}
For even $p\ge 6$ and $q\ge 5$, the tessellation $\{p,q\}$ tolerates combined faults.

\label{Proof:}
For even $p\ge 6$ and $q\ge 5$, let $A$ be the cellular automaton based on the tessellation 
$\{p,q\}$ with symmetric majority voting.
Let $a$ be a cell in $A$, and take $a$ as the origin of the tessellation.
We define the weakened automaton $B$ by having each cell $b$ ignore (1) its parent or parents
(unless $b$ is the origin) and (2) itself (if $q$ is even).
Since $q$ is at least $5$ ($6$ if $q$ is even), the threshold of each vertex is at least $3$ ($4$ if 
$q$ is even).
Thus the reduced threshold of each one-parent vertex is at least $2$, and the reduced threshold of each two-parent vertex is at least $1$.

Say that an edge is {\it special\/} if it is directed from a one-parent parent to a two-parent child.
We associate a flow of $1$ dollar with each special edge, and a flow of $3$ dollars with each other edge.

A one-parent vertex has in-flow at most $3$ (from its parent), since an in-edge to a one-parent vertex cannot be special.
It has out-degree at least $2$ (its reduced threshold).
Of these out-edges, at most one can be special (by Lemma 2.4), so the outflow is at least 
$1+3 = 4$.
Thus every one-parent vertex initially has at least $4-3=1$ dollar.

A two-parent vertex has in-degree $2$ (from its parents), and both of these parents must be 
one-parent vertices (by Lemma 2.5).
Thus the in-flow is at most $1+1 = 2$.
It has out-degree at least $1$ (its reduced threshold), and thus its out-flow is at least $3$
(since an out-edge from a $2$-parent vertex cannot be special).
Thus every two-parent vertex initially has at least $3-2=1$ dollar.

A terminal vertex can have in-degree $1$, in which case the in-flow is $3$, or in-degree $2$, in which case
the in-flow is $1+1=2$.
Thus each terminal vertex  finally has at most
$3$ dollars.
We conclude that $3$ times the number of terminals is at least  the number of non-terminals, and thus that we may take $M=1+3/1=4$ in Theorem 4.2.
\QED

\label{Proposition 4.7:}
For odd $p\ge 5$ and $q\ge 5$, the tessellation $\{p,q\}$ tolerates combined faults.

\label{Proof:}
For odd $p\ge 5$ and $q\ge 5$, let $A$ be the cellular automaton based on the tessellation 
$\{p,q\}$ with symmetric majority voting.
Let $a$ be a cell in $A$, and take $a$ as the origin of the tessellation.
Our definition of the weakened automaton $B$ will be more complicated than that for previous propositions, in that the neighbors of a vertex to be ignored will not be determined merely by the type (cousin, non-cousin one-parent, or two-parent) of that vertex.
Instead, we shall classify vertices as {\it weak\/} or {\it strong}, proceeding inductively, generation by generation, starting with the origin.
The origin, and all non-cousin one-parent vertices, will be strong.
All two-parent vertices will be weak.
Finally, cousin vertices will be strong or weak according as their parent is weak or strong.
We now define the weakened automaton $B$ by having each vertex $b$ ignore (1) its parent or parents, (2) its cousin (if it is a weak cousin vertex) and (3) itself (if $q$ is even).
Since $q$ is at least $5$ ($6$ if $q$ is even), the threshold of each vertex is at least $3$ 
($4$ if $q$ is even).
Each strong vertex ignores at most $1$ neighbor, and thus has reduced threshold at least $2$,
while each weak vertex ignores $2$ neighbors, and thus has reduced threshold at least $1$.

First, we observe that a weak cousin-vertex that is a child of a strong cousin-vertex has a weak cousin-vertex as its cousin.
Suppose, to obtain a contradiction, that the strong cousin-vertex $a$ is the parent of weak cousin-vertex $b$, which has as its cousin the strong cousin-vertex $c$.
Let $d$ be the parent of $c$.
Then $d$ must be weak, since $c$ is strong.
But since $b$ and $c$ are cousins, at least one of their parents, $a$ and $d$, must be a non-cousin one-parent vertex (by Lemma 2.6).
The vertex $a$ is a cousin-vertex, so $d$ must be a non-cousin one-parent vertex.
But this contradicts the fact that $d$ is weak.

Second, we observe that
any cousin edge must be directed out of a strong cousin vertex and into a weak cousin vertex.
To see this, note that by Lemma 2.6, two vertices that are cousins of each other  must have least one non-cousin one-parent parent between them.
This non-cousin one-parent vertex is strong, so at least one of the cousins is weak.
But no cousin edge can be directed out of a weak cousin vertex, since weak cousin vertices ignore their cousin neighbors.
Thus a cousin edge must be directed out of a strong cousin vertex, and into a weak one.

We assign to each cousin-edge a flow of $2$.
We assign to each  edge directed into a two-parent vertex a flow of $4$. 
We assign to each parent-child edge directed into a weak cousin vertex 
having a strong cousin a flow of $6$. 
We assign to each parent-child edge directed into a weak cousin vertex 
having a weak cousin a flow of $8$. 
We assign to each other edge (that is, each edge directed into a non-cousin one-parent vertex or a strong cousin vertex) a flow of $9$.

A non-cousin one-parent vertex has in-degree $1$ (from its parent), and thus has in-flow $9$.
It has out-degree at least $2$ (its reduced threshold).
None of its out-edges can be a cousin-edge (since they are directed from a non-cousin 
one-parent vertex), 
and at most one can be directed into a two-parent vertex (by Lemma 2.4).
Thus its out-flow is at least $4+6 = 10$, so
a non-cousin one-parent vertex has initially at least $10-9=1$ dollar.

A two-parent vertex has in-degree at most $2$ (from its parents), and thus has in-flow at most
$4+4=8$.
It has out-degree at least $1$ (its reduced threshold).
None of its out-edges can be a cousin-edge (since two-parent vertices do not have cousins)
or directed into a weak cousin vertex (since a two-parent vertex, which is weak, cannot have a weak cousin vertex as a child). 
Thus its out-flow is at least $9$, 
so a two-parent vertex has initially at least $9-8=1$ dollar.

A strong cousin vertex has in-degree $1$ (from its parent), and thus has in-flow $9$.
It has out-degree at least $2$ (its reduced threshold).
None of its out-edges can be directed into a two-parent vertex (since, by Lemma 2.5, the parents of a two-parent vertex are non-cousin one-parent vertices).
At most one is a cousin-edge, with an out-flow of $2$, and so at least one of its out-edges must be directed into
a non-cousin one-parent vertex (with an out-flow of $9$) or into a weak cousin vertex.
By the first observation above, such a weak cousin vertex must have a weak cousin vertex as its cousin, contributing an out-flow of $8$ from the strong cousin vertex.
Thus its total out-flow is at least $2+8=10$,
so a strong cousin vertex has initially at least $9-8=1$ dollar.

A weak cousin vertex with a strong cousin has in-degree at most $2$ (from its cousin, its parent, or both).
The in-flow from a cousin-edge is  $2$, and the in-flow from a parent is $6$, so the total in-flow is at most $2+6=8$.
Its out-degree is at least $1$ (its reduced threshold).
None of its out-edges can be a cousin-edge (by the second observation above, weak cousin vertices do not have cousin-edges directed out of them), directed into a two-parent vertex (by Lemma 2.5, the parents of two-parent vertices are non-cousin one-parent vertices), or directed into a weak cousin vertex (a cousin vertex that is the child of a weak cousin vertex is strong).
Thus it must be directed into a non-cousin one-parent vertex or a strong cousin vertex, with an 
out-flow of $9$.
Thus a weak cousin-vertex with a strong cousin has initially at least $1$ dollar.

A weak cousin vertex with a weak cousin has in-degree $1$ (from its parent, since by the second observation above, no edge can be directed to it from its weak cousin).
The in-flow from its parent is $8$. 
Its out-degree is at least $1$ (its reduced threshold).
As in the preceding paragraph, an out-edge must be directed into a non-cousin one-parent vertex or a strong cousin vertex, with an out-flow of $9$.
Thus a weak cousin vertex with a weak cousin has initially at least $1$ dollar.

The analysis above shows that the in-flow of any vertex is at most $9$.
Thus each terminal vertex  finally has at most
$9$ dollars.
We conclude that $9$ times the number of terminals is at least  the number of non-terminals, and thus that we may take $M=1+9/1=10$ in Theorem 4.2.
\QED

Propositions 4.3, 4.4, 4.5, 4.6 and 4.7 together prove the following theorem.

\label{Theorem 4.8:}
The tessellation $\{p,q\}$ tolerates combined faults if $p\ge 5$ and $q\ge 5$,
if $p=4$ and $q\ge 7$, or if $p=3$ and $q\ge 9$.
\sk

\heading{4. Conclusion}

We have classified the regular tessellations of the plane for which binary cellular automata using symmetric majority voting tolerate transient faults or combined faults.
These results can be extended to a much broader class of tessellations: McCann [M2] 
has shown that cellular automata using majority voting and based on ``nice'' graphs
tolerate transient or combined faults if the face and vertex degrees merely satisfy appropriate
upper and lower bounds, together with a technical condition on the parities of face degrees.
(A simple undirected graph is ``nice'' if it is connected, locally-finite, and discretely embeddable in the plane.)

We should also point out that in our  results we have not considered
any transition functions other than those based on majority voting among all neighbors,
which is symmetric under all automorphisms of the underlying graph.
It is known, however, that in other contexts (see Pippenger [P])
asymmetric transition functions are able
to achieve fault tolerance in some situations in which symmetric functions cannot.
\sk

\heading{5. Acknowledgment}

The research reported here was supported
by Grants CCF 0430656 and CCF 0646682 from the National Science Foundation.
\sk

\heading{6.  References}

\refbook B; E. R. Berlekamp, J. H. Conway and R. K. Guy;
Winning Ways for Your Mathematical Plays; Academic Press, 1982, v.~2.

%
\item{[C1]} H. S. M. Coxeter,
``Regular Honeycombs in Hyperbolic Space'', in:
{\it Proceedings of the International Congress of Mathematicians, 1954},
North-Holland Publishing, 1956, v.~III, pp.~155--169
(reprinted in H.~S.~M. Coxeter, {\it The Beauty of Geometry},  Dover Publications, 1999).

\refbook C2; H. S. M. Coxeter and W. O. J. Moser;
Generators and Relations for Discrete Groups;
Springer-Verlag, 1980.

\ref G1; P. G\'{a}cs and J. H. Reif;
``A Simple Three-Dimensional Real-Time Reliable Cellular Array'';
Journal of Computer and System Sciences; 36:2 (1988) 125--147.

\ref G2; M. Gardner;
``Mathematical Games'';
Scientific American; 223:4 (10/1970) 120--123
(reprinted in {\it Wheels, Life, and Other Mathematical Amusements},
W.~H. Freeman and Company, 1983).

\ref M1; M. Marganstern;
``New Tools for Cellular Automata in the Hyperbolic Plane'';
Journal of Universal Computer Science; 6:12 (2000) 1226--1232.

\refbook M2; M. A. McCann;
Memory in Media with Manufacturing Faults;
Ph.~D. Thesis, Department of Computer Science,
Princeton University, September 2007.

\ref M3; M. A. McCann and N. Pippenger;
``Fault Tolerance in Cellular Automata at High Fault Rates'';
Journal of Computer and System Sciences;  74 (2008) 910--918.

\refinbook N1; J. von Neumann;
``Probabilistic Logics and the Synthesis of Reliable Organisms from Unreliable Components'';
in: C.~E. Shannon and J.~McCarthy (Ed's);
Automata Studies; Princeton University Press, 1956, pp.~43--98.

\refbook N2; J. von Neumann (compiled by A.~W. Burks);
 Theory of Self-Reproducing Automata;
 University of Illinois Press, 1966.

\ref P; N. Pippenger;
``Symmetry in Self-Correcting Cellular Automata'';
Journal of Computer and System Sciences; 49:1 (194) 83--95.

\ref S; O. N. Stavskaya and I. I. Pyatetski\u{\i}-Shapiro;
``On Homogeneous Nets of Spontaneously Active Elements'';
Systems Theory Research; 20 (1976) 75--88
 (translation of {\it Problemy Kibernetiki}, 20 (1968) 91--106).

 \ref T1; A. L. Toom;
 ``Nonergodic Multidimensional Systems of Automata'';
 Problems of Information Transmission; 10:3 (1974) 239--246
 (translated from {\it Problemy Peredachi Informatsi\u{\i}}, 10:3 (1974) 70--79).

 \refinbook T2; A. L. Toom;
 ``Stable and Attractive Trajectories in Multicomponent Systems'';
 in: R.~L. Dobrushin and Ya.~A. Sinai (Ed's);
 Multicomponent Random Systems; Marcel Dekker, 1980, pp.~549--575.

 \refinbook U; S. Ulam;
 ``Random Processes and Transformations'';
 in: L.~M. Graves, E.~Hille, P.~Smith and O.~Zariski (Ed's);
 Proceedings of the International Congress of Mathematicians, 1950;
 American Mathematical Society, 1952, v.~2, pp.~264--275.

 \ref W; E. G. Wagner;
 ``On Connecting Modules Together to Form a Modular Computer'';
 IEEE Transactions on Electronic Computers; 15:6 (1966) 864--873.
 
\bye